\documentclass[a4paper,11pt,oneside]{article}
\usepackage[english]{babel}
\usepackage[dvips]{graphicx}
\usepackage{amsthm}
\usepackage{exscale}
\usepackage{amsfonts}
\usepackage{amssymb}
\usepackage[intlimits,sumlimits]{amsmath}
\usepackage{amsmath,cases}
\usepackage[latin1]{inputenc}
\usepackage{mathrsfs}
\usepackage{fancyhdr}

\theoremstyle{definition}

\newtheorem{thm}{Theorem}[section]
\newtheorem{prop}[thm]{Proposition}
\newtheorem{lm}[thm]{Lemma}

\newcommand{\vartot}{\mathrm{d}_{\mathrm{TV}}}
\newcommand{\xb}{\mathbf{x}}
\newcommand{\yb}{\mathbf{y}}
\newcommand{\vb}{\mathbf{v}}
\newcommand{\wb}{\mathbf{w}}
\newcommand{\ub}{\mathbf{u}}
\newcommand{\qb}{\mathbf{q}}
\newcommand{\xib}{\boldsymbol{\xi}}
\newcommand{\psib}{\boldsymbol{\psi}}

\newcommand{\omb}{\boldsymbol{\omega}}
\newcommand{\ud}{\mathrm{d}}

\newcommand{\dxi}{\mathrm{d} \boldsymbol{\xi}}

\newcommand{\unifSo}{u_{S^2}(\ud \boldsymbol{\omega})}

\newcommand{\intethree}{\int_{\mathbb{R}^3}}
\newcommand{\lnorm}{\mid\mid \!}
\newcommand{\rnorm}{\! \mid\mid}

\newcommand{\ind}{1 \! \textrm{l}}
\newcommand{\prob}{\mathscr{P}(\mathbb{R}^3)}
\newcommand{\probTV}{\mathscr{P}_{\mathrm{TV}}(\mathbb{R}^3)}
\newcommand{\borelthree}{\mathscr{B}(\mathbb{R}^3)}
\newcommand{\rone}{\mathbb{R}}
\newcommand{\rthree}{\mathbb{R}^3}

\newcommand{\Lone}{\mathrm{L}^1(\mathbb{R}^3)}
\newcommand{\Cb}{\mathrm{C}_{b}(\mathbb{R}^3)}

\newcommand{\pp}{\textsf{P}}
\newcommand{\pt}{\textsf{P}_t}
\newcommand{\et}{\textsf{E}_t}

\newcommand{\treen}{\mathfrak{t}_n}
\newcommand{\ptreen}{\mathbb{P}(\mathfrak{t}_n)}
\newcommand{\tree}{\mathfrak{t}}
\newcommand{\trees}{\mathfrak{s}}
\newcommand{\treesn}{\mathfrak{s}_{n-1}}
\newcommand{\treenk}{\mathfrak{t}_{n, k}}

\author{\textsc{By Emanuele Dolera and Eugenio Regazzini} \\
\emph{Universit\`a di Modena e Reggio Emilia and Universit\`a di Pavia}}

\title{\textbf{Probabilistic representation for the solution of the homogeneous Boltzmann equation for Maxwellian molecules}}

\date{}

\begin{document}

\maketitle
\begin{abstract}
{\footnotesize Consider the homogeneous Boltzmann equation for Maxwellian molecules. We provide a new representation for its solution in the form of expectation of a random probability measure $\mathcal{M}$. We also prove that the Fourier transform of $\mathcal{M}$ is a conditional characteristic function of a sum of independent random variables, given a suitable $\sigma$-algebra. These facts are then used to prove a CLT for Maxwellian molecules, that is the statement of a necessary and sufficient condition for the weak convergence of the solution of the equation. Such a condition reduces to the finiteness of the second moment of the initial distribution $\mu_0$. As a further application, we give a refinement of some inequalities, due to Elmroth, concerning the evolution of the moments of the solution.}
\end{abstract}


\section{Introduction} \label{sect:intro}

The homogeneous Boltzmann equation in the kinetic theory of dilute gases governs the time evolution of a probability density function $f(\cdot, t)$, which describes, in good approximation, the frequency
$$
\int_{D} f(\vb, t) \ud \vb \approx
\frac{\small{\text{the \ number \ of \ molecules \ with \ velocity \ in} \ D \ \text{at \ time} \ t}}{\small{\text{the \ total \ number \ of \ molecules}}}
$$
for every $D$ in $\borelthree$, the Borel class in $\rthree$. The present treatment is concerned with a \emph{spatially homogeneous gas} composed of a very large number of like particles, moving according to the laws of classical mechanics and colliding in pairs. When there is no outside force the equation is
\begin{eqnarray}
\frac{\partial}{\partial t} f(\mathbf{v}, t) &=& \int_{\mathbb{R}^3}\int_{S^2} [f(\vb_{\ast}, t)
f(\wb_{\ast}, t) \ - \ f(\vb, t) f(\wb, t)] \times \nonumber \\
&\times& B\left(|\wb - \vb|, \frac{\wb - \vb}{|\wb - \vb|} \cdot \omb \right) \unifSo \ud \wb \label{eq:boltzmann}
\end{eqnarray}
with $(\vb, t)$ in $\rthree \times (0, +\infty)$. In the above formula, $u_{S^2}$ stands for the uniform probability measure on the unit sphere $S^2$, embedded in $\rthree$. Moreover, $\vb_{\ast}$ and $\wb_{\ast}$ denote the \emph{post-collisional velocities}, which must obey the conservation of momentum and kinetic energy, that is
$$
\vb + \wb = \ \vb_{\ast} + \wb_{\ast} \ \ \ \ \ \ \ \ \text{and} \ \ \ \ \ \ \ \
|\vb|^2 + |\wb|^2 = \ |\vb_{\ast}|^2 + |\wb_{\ast}|^2 \ .
$$
All such pairs of vectors can be parametrized by unit vectors $\omb$ in $S^2$ and the specific parametrization used throughout this paper -- the so-called $\omb$-representation -- is given by
\begin{equation} \label{eq:velpost}
\begin{array}{lll}
\vb_{\ast} &= \ \vb &+ \ \ [(\wb - \vb) \cdot \omb] \ \omb  \\
\wb_{\ast} &= \ \wb &- \ \ [(\wb - \vb) \cdot \omb] \ \omb
\end{array}
\end{equation}
where $\cdot$ denotes the standard scalar product. The function $B$ is the so-called \emph{collision kernel} and describes the collision process between the molecules at a microscopic level. A complete description of the Boltzmann model can be found in \cite{cerS}. See also the more recent \cite{vil} and references therein.

\emph{All the results in this paper are obtained under two conditions on} $B$. The former is that we are considering \emph{Maxwellian molecules}, meaning that the kernel $B$ does not depend on $|\wb - \vb|$, according to the law
\begin{equation} \label{eq:maxwell}
B\left(|\wb - \vb|, \frac{\wb - \vb}{|\wb - \vb|} \cdot \omb \right) = b\left(\frac{\wb - \vb}{|\wb - \vb|} \cdot \omb \right)\ .
\end{equation}
The so-called \emph{angular collision kernel} $b : (-1, 1) \rightarrow [0, +\infty)$ is a measurable function satisfying the symmetry assumption
\begin{equation} \label{eq:bsymm}
b(x) = b(\sqrt{1 - x^2}) \frac{|x|}{\sqrt{1 - x^2}} = b(-x)
\end{equation}
for all $x$ in $(-1, 1)$, which reflects the indistinguishability of the two colliding particles. The latter condition is usually referred to as the \emph{Grad angular cutoff}, according to which $b$ is required to be integrable. Throughout this paper, the Grad angular cutoff assumption is taken in the form
\begin{equation} \label{eq:cutoff}
\int_{0}^{1} b(x) \ud x = 1 \ .
\end{equation}
Maxwell \cite{max} deduced (\ref{eq:maxwell}) by considering the specific physical interaction originated by a repulsive force proportional to $r^{-5}$, $r$ being the distance of the two colliding particles. In this very particular situation, the kernel $b$ is made explicit, but its analytical form is overcomplicated and possesses a non-integrable singularity.
A few details are contained in Section 5 of Chapter II of \cite{cerS}. On the other hand, the hypotheses of the Grad angular cutoff -- though it allows a wide generality for the function $b$ -- excludes the more meaningful form studied by Maxwell and, for this reason, the fictitious model which is derived is often called \emph{pseudo-Maxwellian}. Nevertheless, the mathematical theory of pseudo-Maxwellian molecules is quite interesting and reveals a very strong link with probability theory, as firstly pointed out by McKean \cite{mck6, mck7}. \\

Now, let us consider some analytical properties of (\ref{eq:boltzmann}), under the assumptions (\ref{eq:maxwell})-(\ref{eq:cutoff}). Given a probability density function $f_0$ on $\rthree$ as initial datum, the resulting Cauchy problem admits a unique solution $f(\cdot, t)$ for each $t \geq 0$, as proved in \cite{mor}. Each $f(\cdot, t)$ is then a probability density function on $\rthree$ for every $t \geq 0$ and the mapping $t \mapsto f(\cdot, t)$ is in $\mathrm{C}([0, +\infty), \Lone) \cap \mathrm{C}^1((0, +\infty), \Lone)$ and $f(\cdot, 0) = f_0(\cdot)$. To allow the initial datum to be as general as possible, the classical framework can be extended to a weak formulation based on probability measures. Maintaining the validity of (\ref{eq:maxwell})-(\ref{eq:cutoff}), the weak form of equation (\ref{eq:boltzmann}) reads
\begin{eqnarray}
\frac{\ud}{\ud t} \intethree \phi(\vb) \mu(\ud \vb, t) &=& \intethree \intethree \int_{S^2} \phi(\vb_{\ast}) b\left(\frac{\wb - \vb}{|\wb - \vb|} \cdot \omb \right) \times \nonumber \\
&\times& \unifSo \mu(\ud \vb, t) \mu(\ud \wb, t) - \intethree \phi(\vb) \mu(\ud \vb, t) \ . \label{eq:wboltzmann}
\end{eqnarray}
The function $\mu$ on $\borelthree \times [0, +\infty)$ is defined to be a \emph{measure solution} when
\begin{enumerate}
\item[i)] the map $t \mapsto \mu(\cdot, t)$ is in  $\mathrm{C}([0, +\infty), \probTV) \cap \mathrm{C}^1((0, +\infty), \probTV)$
\item[ii)] $\mu$ satisfies (\ref{eq:wboltzmann}) for every $t > 0$ and every  test function $\phi$ in $\Cb$, the space of bounded and continuous functions on $\rthree$.
\end{enumerate}
Here, $\probTV$ is the space of all probability measures on $\borelthree$ endowed with the topology generated by the \emph{total variation distance} $\vartot$, which is defined for any pair $(\alpha, \beta)$ of probability measures by
$$
\vartot(\alpha(\cdot); \beta(\cdot)) := \sup_{B \in \borelthree} | \alpha(B) - \beta(B)| \ .
$$
Under hypotheses (\ref{eq:maxwell})-(\ref{eq:cutoff}), if a probability measure $\mu_0$ on $\borelthree$ is given as initial datum, then there exists a unique measure solution $\mu$ such that $\mu(\cdot, 0) = \mu_0(\cdot)$. This extension is discussed in \cite{puto} but, for reader's convenience, some details are included in Subsection \ref{sect:exun} of the present paper. Throughout, the term \emph{solution of} (\ref{eq:boltzmann}) will always mean measure solution, as defined above. \\

This paper proves the validity of a new representation for solutions of (\ref{eq:boltzmann}), in terms of \emph{random probability measures}. The meaning of this locution is specified in Chapter 12 of \cite{ka}. The construction of such a representation has a long story, started with a work of Wild \cite{wil} and then elaborated by McKean \cite{mck6, mck7}, who gave the original Wild formula a probabilistic meaning. A recent development of this kind of research is contained in \cite{blm, blr, cc, ccg0, cgr, dophd, dgr, dore, gr6, gr8}. An outstanding motivation at the basis of these works is the study of the convergence to equilibrium, as $t$ goes to infinity, of solutions of Boltzmann-like equations. In these papers, a fundamental role is played by the so-called \emph{Wild-McKean sum}
\begin{equation} \label{eq:McKean}
\mu(\cdot, t) = \sum_{n = 1}^{+\infty} e^{-t} (1 - e^{-t})^{n-1} \sum_{\treen \in \mathbb{T}(n)} p_n(\treen) \mathcal{Q}_{\treen}[\mu_0](\cdot)
\end{equation}
which will be explained completely in Subsection \ref{sect:Q}. Now, suffice it to recall that every $\mathcal{Q}_{\treen}[\mu_0]$ is a probability measure obtained as a particular $n$-fold convolution of $\mu_0$ with itself, whose definition involves also the kernel $b$. In the one-dimensional case (\emph{Kac's model}), the involvement of $b$ does not alter the nature of $\mathcal{Q}_{\treen}[\mu_0]$ as probability law of a sum of independent random variables. On the contrary, in the three-dimensional case the presence of $b$ in each convolution represents an obstacle for a direct interpretation of the aforementioned type. Since the opportunity of dealing with sums of random variables has proved very fruitful in the study of the Kac model (see \cite{dgr, dore}), then one can try to go round to the aforesaid obstacle in order to recover a new entity interpretable as distribution of a sum. \emph{This is the real motivation for the present work}, whose main results we are going to state.

The first consists in the announced probabilistic representation in terms of a random probability measure, i.e. a random element taking values in the space $\prob$ of all probabilities on $\borelthree$ endowed with the topology of weak convergence.
\begin{thm} \label{thm:main}
\emph{Assume that} (\ref{eq:maxwell})-(\ref{eq:cutoff}) \emph{are in force and that} $\mu$ \emph{is the solution of}
(\ref{eq:boltzmann}) \emph{with initial datum} $\mu_0$. \emph{Then, there exist a measurable space} $(\Omega, \mathscr{F})$, \emph{a family of probability measures} $(\textsf{P}_t)_{t \geq 0}$ \emph{on}
$(\Omega, \mathscr{F})$ \emph{and a random probability measure} $\mathcal{M} : \Omega \rightarrow \prob$ \emph{such that}
\begin{enumerate}
\item[i)] \emph{equality}
\begin{equation} \label{eq:main}
\mu(D, t) = \textsf{E}_t \left[\mathcal{M}(D)\right]
\end{equation}
\emph{holds true for every} $t \geq 0$ \emph{and every} $D \in \borelthree$, \emph{where} $\textsf{E}_t$ \emph{stands for the expectation with respect to} $\textsf{P}_t$.
\item[ii)] $\mathcal{M}$ \emph{depends on} $\mu_0$ \emph{but not on} $b$.
\end{enumerate}
\end{thm}
A precise specification for $(\Omega, \mathscr{F})$ and $(\textsf{P}_t)_{t \geq 0}$ is given in the course of the proof of this theorem, in Subsection \ref{sect:McKean}.

The presentation of the second main result, which contains a representation for $\mathcal{M}$, requires some preliminaries. Define the family $\big{(}(\Omega, \mathscr{F}, \textsf{P}_t)\big{)}_{t \geq 0}$ of probability spaces, mentioned in Theorem \ref{thm:main}, in such a way that they are sufficiently large to support the following random elements. First, a random number $\nu$ taking values in $\mathbb{N} := \{1, 2, \dots\}$. Second, for every $n$ in $\mathbb{N}$, an array $(\pi_{1, n}, \dots, \pi_{n, n})$ of random numbers $\pi_{j, n}$ which take values in $[-1, 1]$. Third, for every $n$ in $\mathbb{N}$, an array $(\mathrm{O}_{1, n}, \dots, \mathrm{O}_{n, n})$ of random matrices $\mathrm{O}_{j, n}$ taking values in $\mathbb{SO}(3)$, the Lie group of orthogonal matrices with positive determinant. Fourth, a sequence $(\mathbf{X}_j)_{j \geq 1}$ of independent and identically distributed (iid) random vectors in $\rthree$. The laws of all the aforesaid random elements will be specified in Subsection \ref{sect:McKean}. The last objects to be defined are elementary mathematical entities of different nature. For each $\ub$ in $S^2$, let $\mathrm{B}(\ub)$ be an element of $\mathbb{SO}(3)$ whose third column coincides with $\ub$. Then, after setting $\mathbf{e}_3 :=\  ^t(0, 0, 1)$, for every $\ub$ in $S^2$, $n$ in $\mathbb{N}$ and $j$ in $\{1, \dots, n\}$, put
\begin{equation} \label{eq:psijn}
\psib_{j, n}(\ub) := \mathrm{B}(\ub) \mathrm{O}_{j, n} \mathbf{e}_3 \ .
\end{equation}
Finally, if $\hat{\zeta}$ denotes the Fourier transform of the probability measure $\zeta$ on $\borelthree$, namely $
\hat{\zeta}(\xib) := \intethree e^{i \xib \cdot \xb} \zeta(\ud \xb)$ for every $\xib$ in $\rthree$, the representation theorem for $\mathcal{M}$ can be stated as follows.
\begin{thm} \label{thm:Mhat}
\emph{Under the same hypotheses of Theorem \ref{thm:main}, one can construct two $\sigma$-algebras} $\mathscr{G}$ \emph{and} $\mathscr{H}$ \emph{such that}:
\begin{enumerate}
\item[i)] $\mathscr{G} \subsetneqq \mathscr{H} \subsetneqq \mathscr{F}$.
\item[ii)] $\mathcal{M}$, $\nu$ \emph{and the arrays} $(\pi_{1, n}, \dots, \pi_{n, n})$ \emph{are} $\mathscr{G}$-\emph{measurable for every} $n$.
\item[iii)] \emph{The arrays} $(\mathrm{O}_{1, n}, \dots, \mathrm{O}_{n, n})$ \emph{are} $\mathscr{H}$-\emph{measurable for every} $n$ \emph{but not} $\mathscr{G}$-\emph{measurable for} $n \geq 2$.
\item[iv)] \emph{The sequence} $(\mathbf{X}_j)_{j \geq 1}$ \emph{is independent of} $\mathscr{H}$.
\item[v)] \emph{For} $\ub$ \emph{in} $S^2$ \emph{and}
\begin{equation} \label{eq:Su}
S(\ub) := \sum_{j = 1}^{\nu} \pi_{j, \nu} \psib_{j, \nu}(\ub) \cdot \mathbf{X}_j
\end{equation}
\emph{there is a version of}
\begin{equation} \label{eq:characteristic}
\textsf{E}_t \left[e^{i \rho S(\ub)} \ | \ \mathscr{H} \right]
\end{equation}
\emph{which, as a function of} $\rho \in \rone$, \emph{is the characteristic function of a (finite) sum of independent random numbers.}
\item[vi)] \emph{Equality}
\begin{equation} \label{eq:Mhat}
\hat{\mathcal{M}}(\xib) = \textsf{E}_t \left[e^{i \rho S(\ub)} \ | \ \mathscr{G} \right] = \textsf{E}_t \left[\textsf{E}_t \left[e^{i \rho S(\ub)} \ | \ \mathscr{H} \right]\ | \ \mathscr{G}\right]
\end{equation}
\emph{holds for every fixed} $\xib \neq \mathbf{0}$, \emph{with} $\rho := |\xib|$, $\ub := \xib/|\xib|$.
\end{enumerate}
\emph{Finally, the terms of} (\ref{eq:Mhat}) \emph{are independent of the choice of matrix} $\mathrm{B}(\ub)$.
\end{thm}
It is worth noting that v)-vi) answer the problem of recovering an interpretation of $\mathcal{M}$ in terms of sums of independent random variables.

The rest of the paper is organized as follows. Section \ref{sect:proof} is concerned with the proof of Theorems \ref{thm:main} and \ref{thm:Mhat}, along with a precise explanation of the probabilistic apparatus introduced in the present section. Sections \ref{sect:CLT}-\ref{sect:elmroth} contain some relevant consequences of these main results. In particular, in Section \ref{sect:CLT} it is shown that the condition
\begin{equation} \label{eq:secondmoment}
\intethree |\xb|^2 \mu_0(\ud \xb) < +\infty
\end{equation}
is necessary and sufficient in order that $\mu(\cdot, t)$ be weakly convergent, as $t$ goes to infinity. Finally, Section \ref{sect:elmroth} includes an improvement of some classical inequalities due to Elmroth on the evolution of the moments of $\mu(\cdot, t)$.

Theorems \ref{thm:main} and \ref{thm:Mhat} can be used also to obtain bounds for the total variation distance between $\mu(\cdot, t)$ and its equilibrium, under additional conditions on $\mu_0$. This problem is one of the most important and mathematically challenging in the entire theory of the Boltzmann equation. A solution in the case of spatially homogeneous Maxwellian molecules, based on representation (\ref{eq:main}), has been presented in \cite{dophd} and will form the subject of future work which, because of its complexity, cannot be incorporated in the present one due to lack of space.

\section{Proof of the main results} \label{sect:proof}

This section, which is the core of the paper, is split into various subsections. The first recalls some tools, due to McKean, which will come in useful in the course of the work, such as the definition of McKean's tree. Subsection \ref{sect:Q} contains the definition of the $Q$ operator, along with some extensions which are needed to explain (\ref{eq:McKean}). Existence and uniqueness is discussed in Subsection \ref{sect:exun}. Then, in Subsection \ref{sect:C}, some formulas related to the $Q$ operator are manipulated to justify the definition of a new operator, $\mathcal{C}$, which sends a pair of probability measures (pms, from now on) onto another pm. Finally, the probabilistic objects introduced in Section \ref{sect:intro} are specified in Subsection \ref{sect:McKean} and combined, successively, to get the proofs of Theorems \ref{thm:main} and \ref{thm:Mhat}, which are provided in Subsections \ref{sect:proof1} and \ref{sect:proof2}, respectively.

\subsection{McKean's trees} \label{sect:trees}

For every $n$ in $\mathbb{N}$, a McKean tree $\treen$ with $n$ leaves is a \emph{binary tree}, whose points are either \emph{nodes} or \emph{leaves}. First, the only tree $\tree_1$ with a single leaf reduces to one point. For $n \geq 2$, every tree $\treen$ possesses exactly $n$ points with no outward connection, the \emph{leaves}, and exactly $n - 1$ points, each generating two ``children'', a \emph{left} and a \emph{right} one, respectively. The top node, the only which is not generated as a child, is called ``root''. From now on, the leaves will be numbered according to a natural left-to-right order. The number of generations separating a leaf from the root is called \emph{depth} and $\delta_j(\treen)$ indicates the depth of the $j$-th leaf of $\treen$. Finally, let $\mathbb{T}(n)$ denote the (finite) set of all McKean's trees with $n$ leaves and $\mathbb{T} := \textsf{X}_{\substack{n \geq 1}} \mathbb{T}(n)$. In this notation, one can define the following operations on trees.
\begin{enumerate}
\item[a)] \emph{Split-up}, when one erases the root of a tree $\treen$ to obtain two trees, a left tree $\treen^l$ and a right tree $\treen^r$, respectively. Henceforth, $n_l$ ($n_r$, respectively) will denote the number of leaves of $\treen^l$ ($\treen^r$, respectively).
\item[b)] \emph{Recombination}, when two trees $\treen$ and $\tree_m$ are combined so that they become the left tree and the right tree, respectively, of a new tree $\treen \oplus \tree_m$ of $n+m$ leaves.
\item[c)] \emph{Germination}, which consists in appending the two-leaved tree $\tree_2$ to a specific leave (say $k$) of $\treen$. Let $\treenk$ indicate the resulting tree.
\end{enumerate}

For a given tree $\treen$ in $\mathbb{T}(n)$, the set of all its germinations constitutes a subset $\mathbb{G}(\treen)$ of $\mathbb{T}(n+1)$, while $\mathbb{P}(\treen)$ denotes the subset of $\mathbb{T}(n-1)$ composed of those trees which can produce $\treen$ by germination.

An important element to understand (\ref{eq:McKean}) is the ``weight'' $p_n(\treen)$ associated with each tree $\treen$. After setting $p_1(\tree_1) := 1$, the definition of $p_n(\treen)$ is given inductively by putting
\begin{equation} \label{eq:pntn}
p_n(\treen) := \frac{1}{n-1} p_{n_l}(\treen^l) p_{n_r}(\treen^r)
\end{equation}
for every $n \geq 2$. Then, one can note straightforwardly that the weights $\{p_n(\treen)\}_{\treen \in \mathbb{T}(n)}$ form a probability distribution on $\mathbb{T}(n)$, in the sense that they all belong to $[0, 1]$ and
$$
\sum_{\treen \in \mathbb{T}(n)} p_n(\treen) = 1
$$
for every $n$ in $\mathbb{N}$.

\subsection{The $Q$ operator and its extension} \label{sect:Q}

Before introducing new definitions, it is worth recalling a few elementary facts about the Boltzmann equation. First, for every $\ub$ in $S^2$, (\ref{eq:bsymm})-(\ref{eq:cutoff}) imply that
\begin{equation}
\int_{S^2} b(\ub \cdot \omb) \unifSo = 1 \ .
\end{equation}
Second, for every $\omb$ in $S^2$, the map $\mathrm{T}_{\omb} : (\vb, \wb) \mapsto (\vb_{\ast}, \wb_{\ast})$, defined in accordance with (\ref{eq:velpost}), is a linear diffeomorphism of $\mathbb{R}^6$ into itself. Linearity follows by inspection, while, for the rest, it is enough to notice that $\mathrm{T}^2_{\omb} \equiv \mathrm{Id}_{\mathbb{R}^6}$, whence $\big{|} Jac \ [\mathrm{T}_{\omb}] \big{|} = 1$. Finally, after a direct verification, $|\wb_{\ast} - \vb_{\ast}| = |\wb - \vb|$ and $(\wb_{\ast} - \vb_{\ast}) \cdot \omb = -(\wb - \vb) \cdot \omb$.

At this stage, for every pair $(p, q)$ of probability density functions (pdfs, from now on) on $\rthree$, define the $Q$ operator according to
\begin{equation} \label{eq:Qop}
Q[p, q](\vb) := \intethree\int_{S^2} p(\vb_{\ast})
q(\wb_{\ast}) b\left(\frac{\wb - \vb}{|\wb - \vb|} \cdot \omb \right) \unifSo \ud \wb \ .
\end{equation}
Its main properties are collected in
\begin{prop} \label{prop:propQ}
\emph{Let} (\ref{eq:bsymm})-(\ref{eq:cutoff}) \emph{be valid. Then,} \emph{for any pair} $(p, q)$ \emph{of pdfs on} $\rthree$, $Q[p, q](\cdot)$ \emph{defines a new pdf on} $\rthree$. \emph{Moreover,} $Q[p, q] = Q[q, p]$ \emph{and, for every bounded and continuous} $\phi$, \emph{one has}
\begin{gather}
\intethree \phi(\vb) Q[p, q](\vb) \ud \vb \nonumber \\
= \intethree \intethree \int_{S^2} \phi(\vb_{\ast}) p(\vb)q(\wb) b\left(\frac{\wb - \vb}{|\wb - \vb|} \cdot \omb \right) \unifSo \ud \vb \ud \wb \ . \label{eq:propQ}
\end{gather}
\emph{Lastly, after setting} $f_R(\vb) := R \vb$, \emph{where} $R$ \emph{is any orthogonal} $3 \times 3$ \emph{matrix}, \emph{equality}
\begin{equation} \label{eq:orthogonalpropertyQ}
Q[p \circ f_{R}^{-1}, q \circ f_{R}^{-1}](\vb) = Q[p, q] \circ f_{R}^{-1}(\vb)
\end{equation}
\emph{holds true for almost every} $\vb$ \emph{in} $\rthree$.
\end{prop}
This statement is classic in the theory of the Boltzmann equation and its proof, here omitted, is contained in any specific text on this subject. See \cite{cerS, vil}.

The analysis of (\ref{eq:Qop}) can be considerably simplified by means of a formula due to Bobylev, concerning the use of the Fourier transform. See \cite{bob88, desv, puto} for details. Since this subject is usually reported starting from a different parametrization of the post-collisional velocities -- a choice that influences the form of the ensuing equations -- the Bobylev formula is stated and proved \emph{ex novo}, using the $\omb$-representation, in
\begin{prop}[Bobylev] \label{prop:bobylev}
\emph{Let assumptions} (\ref{eq:bsymm})-(\ref{eq:cutoff}) \emph{be in force. Then, given any pair} $(p, q)$ \emph{of pdfs on} $\rthree$, \emph{equality}
\begin{equation} \label{eq:bobylev}
\hat{Q}[p, q](\xib) = \int_{S^2} \hat{p}(\xib - (\xib \cdot \omb)\omb) \hat{q}((\xib \cdot \omb)\omb) \ b\left(\frac{\xib}{|\xib|} \cdot \omb \right) \unifSo
\end{equation}
\emph{holds true for every} $\xib \neq \mathbf{0}$.
\end{prop}

\emph{Proof}: Fix any $\xib$ in $\rthree\setminus\{\mathbf{0}\}$ and apply (\ref{eq:propQ}) with $\phi(\vb) = e^{i \vb \cdot \xib}$ to obtain
\begin{equation} \label{eq:stepbobylev}
\hat{Q}[p, q](\xib) = \intethree \intethree \int_{S^2} e^{i \vb_{\ast} \cdot \xib} p(\vb)q(\wb) \ b\left(\frac{\wb - \vb}{|\wb - \vb|} \cdot \omb \right) \unifSo \ud\vb \ud\wb \ .
\end{equation}
Then, recall the definition of $\vb_{\ast}$ to write
\begin{eqnarray}
\exp\{i \vb_{\ast} \cdot \xib\} &=& \exp\{i \vb \cdot \xib\} \cdot \exp\{i [(\wb - \vb) \cdot \omb] \cdot [\omb \cdot \xib]\} \nonumber \\
&=& \exp\{i \vb \cdot \xib\} \cdot \exp\left\{i |\wb - \vb| \cdot |\xib| \cdot \left[\frac{(\wb - \vb)}{|\wb - \vb|} \cdot \omb\right] \cdot \left[\omb \cdot \frac{\xib}{|\xib|}\right]\right\} \nonumber \ .
\end{eqnarray}
At this stage, take into account the integral
$$
\int_{S^2} e^{i \vb_{\ast} \cdot \xib} \ b\left(\frac{\wb - \vb}{|\wb - \vb|} \cdot \omb \right) \unifSo
$$
and change the variable according to $\omb = \ ^tO \boldsymbol{\tau}$, where $O$ is a $3 \times 3$ orthogonal matrix such that $O \frac{\wb - \vb}{|\wb - \vb|} = \frac{\xib}{|\xib|}$ and $O \frac{\xib}{|\xib|} = \frac{\wb - \vb}{|\wb - \vb|}$. From a geometrical point of view, such a matrix corresponds to a reflection, in the plane generated by $\xib$ and $(\wb - \vb)$, around the bisectrix of these very same vectors, provided that they are linearly independent. Of course, for any $\lambda$ in $\rone$, $\{ (\vb, \wb) \in \rone^6 \ | \ (\wb - \vb) = \lambda \xib \}$ is a null set under the measure $p(\vb)q(\wb) \ud\vb \ud\wb$, and so the desired linear independence holds almost everywhere. Since the uniform measure $u_{S^2}$ remains the same after the change $\omb = \ ^tO \boldsymbol{\tau}$, one gets
\begin{gather}
\int_{S^2} e^{i \vb_{\ast} \cdot \xib} \ b\left(\frac{\wb - \vb}{|\wb - \vb|} \cdot \omb \right) \unifSo \nonumber \\
= \int_{S^2}  \ \exp\{i \vb \cdot \xib\} \cdot \exp\{i [(\wb - \vb) \cdot \boldsymbol{\tau}] \cdot [\boldsymbol{\tau} \cdot \xib]\} b\left(\frac{\xib}{|\xib|} \cdot \boldsymbol{\tau} \right) u_{S^2}(\ud \boldsymbol{\tau}) \ . \label{eq:changebobylev}
\end{gather}
Substitution of this expression into (\ref{eq:stepbobylev}) and, then, exchange of the integrals yield
\begin{gather}
\intethree \intethree \exp\{i \vb \cdot \xib\} \cdot \exp\{i [(\wb - \vb) \cdot \boldsymbol{\tau}] \cdot [\boldsymbol{\tau} \cdot \xib]\} p(\vb)q(\wb) \ud\vb \ud\wb \nonumber \\
= \hat{p}\big{(}\xib - (\xib \cdot \boldsymbol{\tau})\boldsymbol{\tau}\big{)} \hat{q}\big{(}(\xib \cdot \boldsymbol{\tau}) \boldsymbol{\tau}\big{)} \nonumber
\end{gather}
and hence (\ref{eq:bobylev}). \ $\Box$ \\

The Bobylev formula is a suitable tool to extend the $Q$ operator from pairs of pdfs $(p, q)$ to pairs of pms. Accordingly, fix two pms $\zeta$ and $\eta$ on $\borelthree$ and, on the same space, choose two sequences $(\zeta_n)_{n\geq 1}$ and $(\eta_n)_{n\geq 1}$ of absolutely continuous pms with densities $(p_n)_{n\geq 1}$ and $(q_n)_{n\geq 1}$, respectively, in such a way that $\zeta_n \rightarrow \zeta$ and $\eta_n \rightarrow \eta$, in the sense of weak convergence of pms. See, e.g., Lemma 9.5.3 in \cite{du}. Finally, set
\begin{equation} \label{eq:extendedQ}
\mathcal{Q}[\zeta, \eta](\ud \vb) := \text{w-lim}_{n \rightarrow \infty} Q[p_n, q_n](\vb) \ud \vb
\end{equation}
which, in view of the next proposition, is given as definition of the desired extension of the $Q$ operator.
\begin{prop} \label{prop:extendedQ}
\emph{The limit in} (\ref{eq:extendedQ}) \emph{exists in the sense of weak convergence of pms and is independent of the choice of the two approximating sequences. The} $\mathcal{Q}$ \emph{operator is continuous in its argument, that is, given two sequences} $(\zeta_n)_{n\geq 1}$ \emph{and} $(\eta_n)_{n\geq 1}$ \emph{(not necessarily absolutely continuous) converging weakly to} $\zeta$ \emph{and} $\eta$, \emph{respectively, one has}
\begin{equation} \label{eq:continuityQ}
\lim_{n \rightarrow \infty} \mathcal{Q}[\zeta_n, \eta_n] = \mathcal{Q}[\zeta, \eta]
\end{equation}
\emph{again in the sense of weak convergence of pms. Moreover, the Bobylev formula extends through the relation}
\begin{equation}
\hat{\mathcal{Q}}[\zeta, \eta](\xib) = \int_{S^2} \hat{\zeta}(\xib - (\xib \cdot \omb)\omb) \hat{\eta}((\xib \cdot \omb)\omb) \ b\left(\frac{\xib}{|\xib|} \cdot \omb \right) \unifSo \ . \label{eq:extendedbobylev}
\end{equation}
\emph{Lastly, after setting} $f_R(\vb) := R \vb$, \emph{where} $R$ \emph{is any orthogonal} $3 \times 3$ \emph{matrix}, \emph{equality}
\begin{equation} \label{eq:orthogonalpropertyQextended}
\mathcal{Q}[\zeta \circ f_{R}^{-1}, \eta \circ f_{R}^{-1}](D) = Q[\zeta, \eta] \circ f_{R}^{-1}(D)
\end{equation}
\emph{holds true for every} $D$ \emph{in} $\borelthree$.
\end{prop}

\emph{Proof}: Start from (\ref{eq:extendedQ}) and apply (\ref{eq:propQ}) with $\phi(\vb) = e^{i \vb_{\ast} \cdot \xib}$, for fixed $\xib$ in $\rthree$, to get
\begin{eqnarray}
\hat{Q}[p_n, q_n](\xib) &=& \intethree \intethree \int_{S^2} e^{i \vb_{\ast} \cdot \xib} p_n(\vb)q_n(\wb) \ b\left(\frac{\wb - \vb}{|\wb - \vb|} \cdot \omb \right) \unifSo \ud\vb \ud\wb \nonumber \\
&=& \intethree \intethree \left[\int_{S^2} e^{i [\omb \cdot \xib] [(\wb - \vb) \cdot \omb]} \ b\left(\frac{\wb - \vb}{|\wb - \vb|} \cdot \omb \right) \unifSo\right] e^{i \vb \cdot \xib} \times \nonumber \\
&\times& p_n(\vb)q_n(\wb)\ud\vb \ud\wb \nonumber \ .
\end{eqnarray}
Now, $\hat{Q}[p_n, q_n](\xib)$ converges pointwise for any $\xib$ when $n$ goes to infinity, from the mapping Theorem 2.7 in \cite{bill2}. The sole point which needs some work is to check that the discontinuities of
$$
(\vb, \wb) \mapsto \int_{S^2} e^{i [\omb \cdot \xib] [(\wb - \vb) \cdot \omb]} \ b\left(\frac{\wb - \vb}{|\wb - \vb|} \cdot \omb \right) \unifSo
$$
form a $\zeta \otimes \eta$ null set $D \subset \mathbb{R}^6$. To this aim, fix a one-dimensional subspace $\ell$ of $\rthree$ in such a way that
\begin{equation} \label{eq:ell}
\zeta \otimes \eta \left(\{(\xb, \yb) \in \mathbb{R}^6 \ | \ \xb - \yb \in \ell_0 \}\right) = 0
\end{equation}
where $\ell_0 := \ell\setminus\{\mathbf{0}\}$. Next, determine an appropriate positive orthonormal basis $\{\mathbf{a}(\ub), \mathbf{b}(\ub), \ub/|\ub| \}$ of $\rthree$ whose components vary with continuity on $\rthree\setminus\ell$ with respect to $\ub := \wb - \vb$, which is possible in view of the so-called ``\emph{hairy ball theorem}''. See \cite{hir} for details. Then, introduce the parametrization
\begin{equation} \label{eq:parametrization}
\omb(\ub) = \sin\varphi \cos\theta \mathbf{a}(\ub) + \sin\varphi \sin\theta \mathbf{b}(\ub) + \cos\varphi \ub/|\ub|
\end{equation}
by which the spherical integral becomes
$$
\frac{1}{4\pi} \int_{0}^{2\pi}\int_{0}^{\pi} e^{i [\omb(\ub) \cdot \xib] |\ub| \cos\varphi} \ b(\cos\varphi) \sin\varphi \ud\varphi \ud\theta \ .
$$
It is now clear that this last integral, as a function of $\ub$, is continuous on $\rthree \setminus \ell_0$ and hence
\begin{equation} \label{eq:Qfourier}
\lim_{n \rightarrow +\infty} \hat{Q}[p_n, q_n](\xib) = \intethree \intethree \int_{S^2} e^{i \vb_{\ast} \cdot \xib} \ b\left(\frac{\wb - \vb}{|\wb - \vb|} \cdot \omb \right) \unifSo \zeta(\ud\vb) \eta(\ud\wb) \ .
\end{equation}
Since this limit, as a function of $\xib$, is clearly continuous, the continuity theorem for characteristic functions can be applied to state that the pm with density $Q[p_n, q_n]$ converges weakly to a pm designated, from now on, by $\mathcal{Q}[\zeta, \eta]$. Finally, observe that the limit is independent of the choice of the approximating sequences $(\zeta_n)_{n\geq 1}$ and $(\eta_n)_{n\geq 1}$ and, therefore, $\mathcal{Q}[\zeta, \eta]$ is well-defined.

To prove (\ref{eq:extendedbobylev}), start from (\ref{eq:Qfourier}) and operate on the spherical integral by applying the argument used by Bobylev to write (\ref{eq:changebobylev}). This way, the conclusion follows exactly as at the end of the proof of Proposition \ref{prop:bobylev}.

The continuity of the $\mathcal{Q}$ operator is an immediate consequence of (\ref{eq:extendedbobylev}): Indeed, it is enough to combine the dominated convergence theorem with the continuity theorem for characteristic functions.

Finally, (\ref{eq:orthogonalpropertyQextended}) can be proved by approximation. Choose two sequences $(\zeta_n)_{n\geq 1}$ and $(\eta_n)_{n\geq 1}$ of absolutely continuous pms, with densities $(p_n)_{n\geq 1}$ and $(q_n)_{n\geq 1}$, respectively, in such a way that $\zeta_n \rightarrow \zeta$ and $\eta_n \rightarrow \eta$ weakly. From (\ref{eq:orthogonalpropertyQ}), $\mathcal{Q}[\zeta_n \circ f_{R}^{-1}, \eta_n \circ f_{R}^{-1}](D) = Q[\zeta_n, \eta_n] \circ f_{R}^{-1}(D)$ for every $D$ in $\borelthree$. At this point, first observe that if a sequence $(\pi_n)_{n \geq 1}$ of pms converges weakly to the pm $\pi$, then also $(\pi_n \circ h^{-1})_{n \geq 1}$ converges weakly to $\pi \circ h^{-1}$, for every continuous function $h: \rthree \rightarrow \rthree$. Thus, to get (\ref{eq:orthogonalpropertyQextended}), exploit the weak continuity and take the limit in $n$ of both sides. $\Box$ \\

It remains to define the pms $\mathcal{Q}_{\treen}[\mu_0]$, which completes the description of (\ref{eq:McKean}). After fixing the pm $\mu_0$ on $\borelthree$, put
\begin{equation} \label{eq:Qtn}
\begin{array}{lll}
\mathcal{Q}_{\tree_1}[\mu_0] &:= \mu_0 & {} \\
\mathcal{Q}_{\treen}[\mu_0] &:= \mathcal{Q}\left[\mathcal{Q}_{\treen^l}[\mu_0], \mathcal{Q}_{\treen^r}[\mu_0]\right] & \ \ \ \ \text{for} \ n \geq 2 \ .
\end{array}
\end{equation}

\subsection{Existence and uniqueness for the Cauchy problem} \label{sect:exun}

For the existence, it is enough to show that the right-hand side (RHS, in short) of (\ref{eq:McKean}) gives a solution of (\ref{eq:boltzmann}) when hypotheses (\ref{eq:maxwell})-(\ref{eq:cutoff}) are in force. Of course, when $t = 0$, the series reduces to $\mu_0$ and so the initial condition is fulfilled. The regularity of the map
$$
t \rightarrow \sum_{n = 1}^{+\infty} e^{-t} (1 - e^{-t})^{n-1} \sum_{\treen \in \mathbb{T}(n)} p_n(\treen) \mathcal{Q}_{\treen}[\mu_0] =: \mathcal{W}_t[\mu_0]
$$
follows from the fact that $\vartot(\mathcal{W}_{t_1}[\mu_0]; \mathcal{W}_{t_2}[\mu_0]) \leq c_0(t_1 t_2) |t_1 - t_2|$ for every $t_1$ and $t_2$ in $[0, +\infty)$ and
\begin{eqnarray}
\sup_{B \in \borelthree} \Big{|} &\sum_{n = 2}^{+\infty}& (n-1) e^{-t_1} (1 - e^{-t_1})^{n-2} \sum_{\treen \in \mathbb{T}(n)} p_n(\treen) \mathcal{Q}_{\treen}[\mu_0](B) \nonumber \\
- &\sum_{n = 2}^{+\infty}& (n-1) e^{-t_2} (1 - e^{-t_2})^{n-2} \sum_{\treen \in \mathbb{T}(n)} p_n(\treen) \mathcal{Q}_{\treen}[\mu_0](B) \Big{|} \nonumber \\
&\leq& C_1(t_1, t_2) |t_1 - t_2| \nonumber
\end{eqnarray}
for every $t_1$ and $t_2$ in $(0, +\infty)$. Here, $c_0$ and $c_1$ are positive constants for which $\limsup_{t_2 \rightarrow t_1} c_i(t_1, t_2) < +\infty$, for $i = 0, 1$ and every $t_1$ and $t_2$. The expression of these constants can be essentially obtained by the elementary relations $|e^{-t_1} - e^{-t_2}| \leq |t_1 - t_2|$ and $|a^n - b^n| \leq (n-1) \max\{a^{n-1}; b^{n-1}\} |a - b|$, valid for positive $t_1$, $t_2$, $a$ and $b$. Finally, to show that the RHS in (\ref{eq:McKean}) solves (\ref{eq:wboltzmann}) for every $\phi$ in $\mathrm{C}_b$, define, for every $n$ in $\mathbb{N}$,
$$
\begin{array}{lll}
\mathcal{Q}_1[\mu_0] &:= \mathcal{Q}_{\tree_1}[\mu_0] & {} \\
\mathcal{Q}_n[\mu_0] &:= \sum_{\treen \in \mathbb{T}(n)} p_n(\treen) \mathcal{Q}_{\treen}[\mu_0] & \ \ \ \ \text{for} \ n \geq 2
\end{array}
$$
and note that
\begin{equation} \label{eq:lala}
\mathcal{Q}_n[\mu_0] = \frac{1}{n-1} \sum_{k = 1}^{n-1} \mathcal{Q}[\mathcal{Q}_k[\mu_0], \mathcal{Q}_{n-k}[\mu_0]]
\end{equation}
holds for every $n \geq 2$. Thus, differentiation term by term gives
\begin{eqnarray}
&{}& \frac{\ud}{\ud t} \intethree \phi(\xb) \mathcal{W}_t[\mu_0](\ud \xb) \nonumber \\
&=& - \intethree \phi(\xb) \mathcal{W}_t[\mu_0](\ud \xb) + \sum_{n = 2}^{+\infty} (n-1) e^{-2t} (1 - e^{-t})^{n-2} \intethree \phi(\xb) \mathcal{Q}_n[\mu_0](\ud \xb) \nonumber
\end{eqnarray}
and, from (\ref{eq:lala}) and the properties of $\mathcal{Q}$,
\begin{gather}
\sum_{n = 2}^{+\infty} (n-1) e^{-2t} (1 - e^{-t})^{n-2} \intethree \phi(\xb) \mathcal{Q}_n[\mu_0](\ud \xb) \nonumber \\
= \sum_{n = 2}^{+\infty} e^{-2t} (1 - e^{-t})^{n-2} \sum_{k = 1}^{n-1} \intethree \phi(\xb) \mathcal{Q}[\mathcal{Q}_k[\mu_0], \mathcal{Q}_{n-k}[\mu_0]](\ud \xb) \nonumber \\
= \sum_{n = 2}^{+\infty} e^{-2t} (1 - e^{-t})^{n-2} \times \nonumber \\
\times  \sum_{k = 1}^{n-1} \intethree \intethree \int_{S^2} \phi(\vb_{\ast}) b\left(\frac{\wb - \vb}{|\wb - \vb|} \cdot \omb \right) \unifSo \mathcal{Q}_k[\mu_0](\ud \vb) \mathcal{Q}_{n-k}[\mu_0](\ud \wb) \nonumber \\
= \sum_{n = 1}^{+\infty} e^{-t} (1 - e^{-t})^{n-1} \sum_{m = 1}^{+\infty} e^{-t} (1 - e^{-t})^{m-1} \times \nonumber \\
\times \intethree \intethree \int_{S^2} \phi(\vb_{\ast}) b\left(\frac{\wb - \vb}{|\wb - \vb|} \cdot \omb \right) \unifSo \mathcal{Q}_n[\mu_0](\ud \vb) \mathcal{Q}_m[\mu_0](\ud \wb) \nonumber \\
= \intethree \intethree \int_{S^2} \phi(\vb_{\ast}) b\left(\frac{\wb - \vb}{|\wb - \vb|} \cdot \omb \right) \unifSo \mathcal{W}_t[\mu_0](\ud \vb) \mathcal{W}_t[\mu_0](\ud \wb) \nonumber
\end{gather}
which proves the claim.

For the uniqueness, let $\mu^{(1)}(\cdot, t)$ and $\mu^{(2)}(\cdot, t)$ be solutions with respect to the same initial datum $\mu_0$. Put
$$
\Delta_t := \sup_{\substack{\phi \in \mathrm{C}_0(\rthree) \\ \lnorm \ \phi \ \rnorm_{\infty} \ \leq 1}} \Big{|} \intethree \phi(\vb) \mu^{(1)}(\ud \vb, t) - \intethree \phi(\vb) \mu^{(2)}(\ud \vb, t)\Big{|}
$$
where $\lnorm \phi \rnorm_{\infty} := \sup_{\vb \in \rthree} |\phi(\vb)|$ and $\mathrm{C}_0(\rthree)$ is the set of all continuous functions on $\rthree$ which vanish at infinity. In view of the Riesz representation theorem, $\Delta_t = \vartot(\mu^{(1)}(\cdot, t); \mu^{(2)}(\cdot, t))$. See Theorem 7.17 in \cite{fol}. Then, define $F_{\phi}(\vb, \wb) := \int_{S^2} \phi(\vb_{\ast}) b\left(\frac{\wb - \vb}{|\wb - \vb|} \cdot \omb \right) \unifSo$ for every $\phi$ in $\mathrm{C}_0(\rthree)$ with $\lnorm \phi \rnorm_{\infty} \ \leq 1$, so that
\begin{eqnarray}
&{}& \Big{|} \intethree \intethree F_{\phi}(\vb, \wb) [\mu^{(1)}(\ud \vb, t)\mu^{(1)}(\ud \wb, t) - \mu^{(2)}(\ud \vb, t)\mu^{(2)}(\ud \wb, t)] \Big{|} \nonumber \\
&\leq& \Big{|} \intethree \intethree F_{\phi}(\vb, \wb) [\mu^{(1)}(\ud \vb, t)\mu^{(1)}(\ud \wb, t) - \mu^{(1)}(\ud \vb, t)\mu^{(2)}(\ud \wb, t)] \Big{|} \nonumber \\
&+& \Big{|} \intethree \intethree F_{\phi}(\vb, \wb) [\mu^{(1)}(\ud \vb, t)\mu^{(2)}(\ud \wb, t) - \mu^{(2)}(\ud \vb, t)\mu^{(2)}(\ud \wb, t)] \Big{|} \leq 2\Delta_t \nonumber
\end{eqnarray}
in view of Lemma 9.1.1 of \cite{stro} and the fact that $\vb \mapsto \intethree F_{\phi}(\vb, \wb) \mu^{(i)}(\ud \wb, t)$ is bounded for $i = 1, 2$. Integration of (\ref{eq:wboltzmann}) with respect to time combined with the last inequality, gives $\Delta_t \leq 3 \int_{0}^{t} \Delta_s \ud s$. At this stage, the Gronwall lemma entails $\Delta_t = 0$ for all $t \geq 0$, which is tantamount to asserting that $\mu^{(1)}(\cdot, t) = \mu^{(2)}(\cdot, t)$. $\Box$

We conclude this subsection by pointing out a useful consequence of the Wild-McKean sum (\ref{eq:McKean}) combined with (\ref{eq:orthogonalpropertyQextended}) and (\ref{eq:Qtn}).
\begin{lm} \label{lm:orthogonal}
\emph{Assume that} (\ref{eq:maxwell})-(\ref{eq:cutoff}) \emph{are in force. Let} $\mu_0$ \emph{be any pm on} $\borelthree$ \emph{and let} $\mu(\cdot, t)$ \emph{be the solution of} (\ref{eq:boltzmann}) \emph{with initial datum} $\mu_0$. \emph{Then, given any orthogonal} $3 \times 3$ \emph{matrix} $R$ \emph{and the relative map} $\mathcal{R}: \xb \mapsto R\xb$, \emph{the solution of} (\ref{eq:boltzmann}) \emph{with initial datum} $\mu_0 \circ \mathcal{R}^{-1}$ \emph{coincides with} $\mu(\cdot, t) \circ \mathcal{R}^{-1}$.
\end{lm}

\subsection{The $\mathcal{C}$ operator and its extension} \label{sect:C}

The extended Bobylev formula (\ref{eq:extendedbobylev}) is now re-examined to derive the definition of a new convolution between pms, say $\zeta$ and $\eta$, on $\borelthree$. Start from that formula and rewrite the integral by changing the variables as follows: First, for every unit vector $\ub$ in $S^2$, fix two vectors $\mathbf{a}(\ub)$ and $\mathbf{b}(\ub)$ in such a way that $\{\mathbf{a}(\ub), \mathbf{b}(\ub), \ub\}$ is a positive orthonormal basis of $\rthree$, regardless of the regularity of the mappings $\ub \mapsto \mathbf{a}(\ub)$ and $\ub \mapsto \mathbf{b}(\ub)$. Second, for any fixed $\xib \neq \mathbf{0}$, choose the parametrization given by
$$
\omb = \sin\varphi \cos\theta \mathbf{a}(\ub) + \sin\varphi \sin\theta \mathbf{b}(\ub) + \cos\varphi \ub
$$
where $\ub := \xib/|\xib|$ and $(\varphi, \theta)$ belongs to $[0, \pi] \times (0, 2\pi)$, which becomes the new domain of integration. Then, observe that $(\xib \cdot \omb) = \rho \cos \varphi$, with $\rho := |\xib|$, and that the uniform pm $\unifSo$ transforms into $(1/4\pi) \sin\varphi \ud \varphi \ud \theta$. An elementary computation leads to
\begin{equation} \label{eq:parametrizedbobylev}
\hat{\mathcal{Q}}[\zeta, \eta](\xib) = \int_{0}^{\pi}\int_{0}^{2\pi} \hat{\zeta}(\rho \cos\varphi \boldsymbol{\psi}^l) \hat{\eta}(\rho \sin\varphi \boldsymbol{\psi}^r) u_{(0, 2\pi)}(\ud \theta) \beta(\ud \varphi)
\end{equation}
for every $\xib \neq \mathbf{0}$, where $u_{(0, 2\pi)}$ stands for the continuous uniform pm on $(0, 2\pi)$, $\beta$ is the pm on $[0, \pi]$ defined by
\begin{equation} \label{eq:beta}
\beta(\ud \varphi) := \frac{1}{2} b(\cos\varphi) \sin\varphi \ud \varphi
\end{equation}
while $\boldsymbol{\psi}^l$ and $\boldsymbol{\psi}^r$ are abbreviations for the quantities
\begin{equation} \label{eq:psi}
\begin{array}{lll}
\psib^l(\varphi, \theta, \ub) &:= &{}\cos\theta \sin\varphi \mathbf{a}(\ub) + \sin\theta \sin\varphi \mathbf{b}(\ub) + \cos\varphi \ub \\
\psib^r(\varphi, \theta, \ub) &:= &-\cos\theta \cos\varphi \mathbf{a}(\ub) - \sin\theta \cos\varphi \mathbf{b}(\ub) + \sin\varphi \ub \ .
\end{array}
\end{equation}
Observe that the realizations of $\boldsymbol{\psi}^l$ and $\boldsymbol{\psi}^r$ depend crucially on the choice of the basis $\{\mathbf{a}(\ub), \mathbf{b}(\ub), \ub\}$. Now, some distinguishing properties of the function
\begin{equation} \label{eq:internalintegralbobylev}
\mathrm{I}(\xib, \varphi) := \left\{ \begin{array}{ll}
\int_{0}^{2\pi} \hat{\zeta}(\rho \cos\varphi \boldsymbol{\psi}^l) \hat{\eta}(\rho \sin\varphi \boldsymbol{\psi}^r) u_{(0, 2\pi)}(\ud \theta) & \text{if} \ \xib \neq \mathbf{0} \\
1 & \text{if} \ \xib = \mathbf{0} \ ,
\end{array} \right.
\end{equation}
which appears in (\ref{eq:parametrizedbobylev}), are highlighted by
\begin{prop} \label{prop:C}
\emph{For every fixed pairs of pms} $(\zeta, \eta)$, $\mathrm{I}$ \emph{is independent of the choice of the basis} $\{\mathbf{a}(\ub), \mathbf{b}(\ub), \ub\}$, \emph{for every} $\ub$ \emph{in} $S^2$, \emph{and turns out to be a measurable function of} $(\xib, \varphi)$ \emph{in} $\rthree \times [0, \pi]$. \emph{For every fixed} $\varphi$ \emph{in} $[0, \pi]$,
$I(\cdot, \varphi)$ \emph{is the Fourier transform of a pm, say} $\mathcal{C}[\zeta, \eta; \varphi]$, \emph{on} $\borelthree$, \emph{that is}
\begin{equation} \label{eq:bobylevforC}
\mathrm{I}(\xib, \varphi) = \hat{\mathcal{C}}[\zeta, \eta; \varphi](\xib)
\end{equation}
\emph{for every} $\xib$ \emph{in} $\rthree$. \emph{Moreover, the mapping} $\varphi \mapsto \mathcal{C}[\zeta, \eta; \varphi]$ \emph{is a random pm and}
\begin{equation} \label{eq:QCbis}
\mathcal{Q}[\zeta, \eta](D) = \int_{0}^{\pi} \mathcal{C}[\zeta, \eta; \varphi](D) \beta(\ud \varphi)
\end{equation}
\emph{is valid for every} $D$ \emph{in} $\borelthree$.
\end{prop}

\emph{Proof}: To prove the first claim, fix $\xib \neq \mathbf{0}$ and let $\{\mathbf{a}(\ub), \mathbf{b}(\ub), \ub \}$ and $\{\mathbf{a}^{'}(\ub), \mathbf{b}^{'}(\ub), \ub \}$ be two distinct positive bases. Since there exists some $\theta^{\ast}$ in $[0, 2\pi)$ such that
\begin{eqnarray}
\mathbf{a}^{'} &=& \cos\theta^{\ast} \mathbf{a} - \sin\theta^{\ast} \mathbf{b} \nonumber \\
\mathbf{b}^{'} &=& \sin\theta^{\ast} \mathbf{a} + \cos\theta^{\ast} \mathbf{b} \ , \nonumber
\end{eqnarray}
write
$$
\begin{array}{lll}
\boldsymbol{\psi}^l(\varphi, \theta, \ub) &:= &{}\cos\theta \sin\varphi \mathbf{a}^{'}(\ub) + \sin\theta \sin\varphi \mathbf{b}^{'}(\ub) + \cos\varphi \ub \\
{} &= &{}\cos(\theta - \theta^{\ast}) \sin\varphi \mathbf{a}(\ub) + \sin(\theta - \theta^{\ast}) \sin\varphi \mathbf{b}(\ub) + \cos\varphi \ub \\
\boldsymbol{\psi}^r(\varphi, \theta, \ub) &:= &-\cos\theta \cos\varphi \mathbf{a}^{'}(\ub) - \sin\theta \cos\varphi \mathbf{b}^{'}(\ub) + \sin\varphi \ub \\
{} &= &{}-\cos(\theta - \theta^{\ast}) \cos\varphi \mathbf{a}(\ub) - \sin(\theta - \theta^{\ast}) \cos\varphi \mathbf{b}(\ub) + \sin\varphi \ub
\end{array}
$$
and substitute them into the integral in (\ref{eq:internalintegralbobylev}). An obvious change of variable leads to the desired conclusion.

To proceed to the other points, it is first proved that $\xib \mapsto \mathrm{I}(\xib, \varphi)$ is continuous. Start by taking a sequence $(\xib_n)_{n \geq 1}$ converging to $\xib^{\ast}$ and put $\rho_n := |\xib_n|$. First, assume that $\xib^{\ast} = \mathbf{0}$ and, to avoid trivialities, $\rho_n \neq 0$ for every $n$ in $\mathbb{N}$. Then, $\rho_n$ goes to zero along with $|\rho_n \cos \varphi \boldsymbol{\psi}^l(\theta, \varphi, \ub_n)|$ and $|\rho_n \sin \varphi \boldsymbol{\psi}^r(\theta, \varphi, \ub_n)|$, where $\ub_n := \xib_n/\rho_n$, for every $\varphi$ in $[0, \pi]$ and $\theta$ in $(0, 2\pi)$. An application of the dominated convergence theorem shows that $\mathrm{I}(\xib_n, \varphi)$ converges to one, for every $\varphi$ in $[0, \pi]$. If $\xib^{\ast} \neq \mathbf{0}$, observe that $\rho_n$ converges to $\rho^{\ast} := |\xib^{\ast}|$, as well as $\ub_n$ converges to $\ub^{\ast} := \xib^{\ast}/\rho^{\ast}$. Fix a small open neighborhood $\Omega(\ub^{\ast}) \subset S^2$ of $\ub^{\ast}$ (with respect to the standard topology of $S^2$) in such a way that $S^2 \setminus \overline{\Omega(\ub^{\ast})}$ contains at least two antipodal points. In view of the independence of $\mathrm{I}(\xib, \varphi)$ of the basis $\{\mathbf{a}(\ub), \mathbf{b}(\ub), \ub \}$, fix a distinguished basis in such a way that the restrictions of $\ub \mapsto \mathbf{a}(\ub)$ and $\ub \mapsto \mathbf{b}(\ub)$ to $\Omega(\ub^{\ast})$ vary with continuity. As a consequence, $\boldsymbol{\psi}^l(\varphi, \theta, \ub_n)$ converges to $\boldsymbol{\psi}^l(\varphi, \theta, \ub^{\ast})$ and $\boldsymbol{\psi}^r(\varphi, \theta, \ub_n)$ converges to $\boldsymbol{\psi}^r(\varphi, \theta, \ub^{\ast})$, for every $\varphi$ in $[0, \pi]$ and $\theta$ in $(0, 2\pi)$. At this stage, the convergence of $\mathrm{I}(\xib_n, \varphi)$ to $\mathrm{I}(\xib^{\ast}, \varphi)$, for every $\varphi$ in $[0, \pi]$, follows once again by an application of the dominated convergence theorem.

As for the measurability of $(\xib, \varphi) \mapsto \mathrm{I}(\xib, \varphi)$, invoke Proposition 9 in Section 9.3 of \cite{frgr}. In view of the continuity of $\xib \mapsto \mathrm{I}(\xib, \varphi)$ for every $\varphi$ in $[0, \pi]$, it suffices to prove that $\varphi \mapsto \mathrm{I}(\xib, \varphi)$ is measurable for each fixed $\xib$. This claim follows from the continuity of $\varphi \mapsto \mathrm{I}(\xib, \varphi)$, that can be verified by observing the explicit dependence on $\varphi$ in (\ref{eq:psi}) and (\ref{eq:internalintegralbobylev}), regardless of the choice of the basis $\{\mathbf{a}(\ub), \mathbf{b}(\ub), \ub \}$.

To show that $\xib \mapsto \mathrm{I}(\xib, \varphi)$ is a characteristic function for every $\varphi$ in $[0, \pi]$, one can resort to the multivariate version of the Bochner characterization. See Exercise 3.1.9 in \cite{stro}. The only point that requires some care is, at this stage, positivity. If this property were not in force, one could find a positive integer $N$, two $N$-vectors $(\omega_1, \dots, \omega_N)$ and $(\xib_1, \dots, \xib_N)$ in $\mathbb{C}^N$ and $(\rthree)^N$, respectively, and some $\varphi^{\ast}$ in $[0, \pi]$ in such a way that
$$
\sum_{j = 1}^{N} \sum_{k = 1}^{N} \omega_j \overline{\omega}_k \mathrm{I}(\xib_j - \xib_k, \varphi^{\ast}) < 0 \ .
$$
Hence, by continuity of $\varphi \mapsto \mathrm{I}(\xib, \varphi)$, for each fixed $\xib$, there exists an open interval $J$ in $[0, \pi]$ containing $\varphi^{\ast}$ such that
$$
\varphi \mapsto \sum_{j = 1}^{N} \sum_{k = 1}^{N} \omega_j \overline{\omega}_k \mathrm{I}(\xib_j - \xib_k, \varphi)
$$
is strictly negative on $J$. Now, choose a specific Maxwellian kernel $b_{\ast}$ for which the resulting pm in (\ref{eq:beta}), say $\beta_{\ast}$, is supported by $\overline{J}$. By construction,
\begin{equation} \label{eq:bochner1}
\int_{0}^{2\pi} \sum_{j = 1}^{N} \sum_{k = 1}^{N} \omega_j \overline{\omega}_k \mathrm{I}(\xib_j - \xib_k, \varphi) \beta_{\ast}(\ud \varphi)
\end{equation}
is a strictly negative number, a fact which immediately leads to a contradiction. Indeed, define $\mathcal{Q}_{\ast}[\zeta, \eta]$ to be the pm with Fourier transform given by (\ref{eq:extendedbobylev}) with $b_{\ast}$ in place of $b$, and note that, in view of (\ref{eq:parametrizedbobylev}), (\ref{eq:bochner1}) is equal to
$$
\sum_{j = 1}^{N} \sum_{k = 1}^{N} \omega_j \overline{\omega}_k \hat{\mathcal{Q}}_{\ast}[\zeta, \eta](\xib_j - \xib_k)
$$
a quantity that must be non negative, from the Bochner criterion.

To prove the measurability of $\varphi \mapsto \mathcal{C}[\zeta, \eta; \varphi]$ it is enough to verify that each map $\varphi \mapsto \mathcal{C}[\zeta, \eta; \varphi](K)$ is $\mathscr{B}([0, \pi])/\mathscr{B}([0, 1])$-measurable, for every rectangle $K = \textsf{X}_{\substack{i=1}}^{3} (-\infty, x_i]$. See, for example, Lemma 1.40 of \cite{ka}. To this end, note that the Fubini theorem can be applied to show that
\begin{eqnarray}
&{}& (\varphi, \mathbf{b}) \mapsto \textsf{G}(\varphi, \mathbf{b}) \nonumber \\
&:=& \lim_{\mathbf{a} \rightarrow -\boldsymbol{\infty}} \lim_{c \rightarrow +\infty} \frac{1}{(2\pi)^3} \int_{-c}^{c} \int_{-c}^{c} \int_{-c}^{c} \left[\prod_{m=1}^{3} \frac{e^{-i\xi_m a_m} - e^{-i\xi_m b_m}}{i\xi_m} \right] \hat{\mathcal{C}}[\zeta, \eta; \varphi](\xib) \dxi  \nonumber
\end{eqnarray}
is $\mathscr{B}(\rthree \times [0, \pi])/\mathscr{B}([0, 1])$-measurable. See Section 8.5 of \cite{chte}. At this stage, to complete the argument it suffices to note that
$$
\mathcal{C}[\zeta, \eta; \varphi]\big{(}\textsf{X}_{\substack{i=1}}^{3} (-\infty, x_i]\big{)} = \lim_{\mathbf{b} \downarrow \xb} \textsf{G}(\varphi, \mathbf{b}) \ .
$$

Finally, equality (\ref{eq:QCbis}) can be derived from (\ref{eq:parametrizedbobylev}) and (\ref{eq:bobylevforC}), through the following lemma, which completes the proof of the proposition.

\begin{lm} \label{lm:mixturechacteristics}
\emph{Let} $(S, \mathscr{S}, \lambda)$ \emph{be a probability space}, $\pi$ \emph{be a pm on} $\mathscr{B}(\rone^d)$ \emph{and} $\eta : S \times \mathscr{B}(\rone^d) \rightarrow [0, +\infty)$ \emph{be a probability kernel with} $\eta_s(\cdot) := \eta(s, \cdot)$ \emph{for every} $s$ \emph{in} $S$. \emph{If} $\hat{\pi}(\xib) = \int_{S} \hat{\eta}_s(\xib) \lambda(\ud s)$ \emph{for every} $\xib$ \emph{in} $\rone^d$, \emph{then} $\pi(D) = \int_{S} \eta_s(D) \lambda(\ud s)$ \emph{holds true for every} $D$ \emph{in} $\mathscr{B}(\rone^d)$.
\end{lm}

\emph{Proof}: For $\epsilon > 0$, set
$$
\begin{array}{l}
\hat{\pi}^{\epsilon}(\xib) := \hat{\pi}(\xib) \exp\{- \epsilon |\xib|^2/2 \} \\
\hat{\eta}_{s}^{\epsilon}(\xib) := \hat{\eta}_s(\xib) \exp\{- \epsilon |\xib|^2/2 \}
\end{array}
$$
for every $s$ in $S$. Then, equality $\hat{\pi}^{\epsilon}(\xib) = \int_{S} \hat{\eta}^{\epsilon}_s(\xib) \lambda(\ud s)$ is still valid for every $\xib$ in $\rone^d$, which entails $\pi^{\epsilon}(D) = \int_{S} \eta^{\epsilon}_{s}(D) \lambda(\ud s)$ for every $D$ in $\mathscr{B}(\rone^d)$, from the classical inversion theorem for the Fourier transform. Now, $\pi^{\epsilon}$ and $\eta^{\epsilon}_{s}$ converge weakly to $\pi$ and $\eta_{s}$, respectively, as $\epsilon \downarrow 0$. See, for example, Lemma 9.5.3 in \cite{du}. Moreover,
$$
\liminf_{\epsilon \downarrow 0} \pi^{\epsilon}(A) \geq \int_{S} \liminf_{\epsilon \downarrow 0} \eta^{\epsilon}_{s}(A) \lambda(\ud s) \geq \int_{S} \eta_{s}(A) \lambda(\ud s)
$$
where the former inequality follows from the Fatou lemma, while the latter holds for every open subset $A$ of $\rone^d$, in view of Theorem 2.1 (iv) of \cite{bill2} applied to $\eta^{\epsilon}_{s}$. A further application of this very same theorem to $\pi^{\epsilon}$ gives the weak convergence of $\pi^{\epsilon}$ to $\int_{S} \eta_s \lambda(\ud s)$ and the thesis follows from the uniqueness of the weak limit. \ $\Box$

\vspace{5mm}

In view of the last proposition, given any $\varphi$ in $[0, \pi]$, $\mathcal{C}[\zeta, \eta; \varphi]$ can be seen as an operator which sends the pair $(\zeta, \eta)$ onto a new pm. Then, it is natural to mimic the procedure of iteration for the $\mathcal{Q}$ operator, defined in Subsection \ref{sect:Q}, to get an analogous result for the $\mathcal{C}$ operator. After fixing the pm $\mu_0$ on $\borelthree$, set $\mathcal{C}_{\mathfrak{t}_1}[\mu_0; \emptyset] := \mu_0$ and, for every tree $\treen$ with $n \geq 2$ and any vector $\boldsymbol{\varphi} = (\varphi_1, \dots, \varphi_{n-1})$ in $[0, \pi]^{n-1}$, introduce the symbols
$$
\begin{array}{ll}
\boldsymbol{\varphi}^l &:= (\varphi_1, \dots, \varphi_{n_l-1}) \\
\boldsymbol{\varphi}^r &:= (\varphi_{n_l}, \dots, \varphi_{n-2})
\end{array}
$$
to define
\begin{equation} \label{eq:Cmisureiterate}
\mathcal{C}_{\treen}[\mu_0; \boldsymbol{\varphi}] := \mathcal{C}\left[\mathcal{C}_{\treen^l}[\mu_0; \boldsymbol{\varphi}^l], \mathcal{C}_{\treen^r}[\mu_0; \boldsymbol{\varphi}^r]; \varphi_{n-1}\right]
\end{equation}
with the proviso that $\boldsymbol{\varphi}^l$ ($\boldsymbol{\varphi}^r$, respectively) is void when $n_l$ ($n - n_l$, respectively) is equal to one. These definitions, combined with Proposition \ref{prop:C}, lead to
\begin{prop} \label{prop:Ctn}
\emph{For every tree} $\treen$ \emph{in} $\mathbb{T}$, \emph{the mapping} $\boldsymbol{\varphi} \mapsto \mathcal{C}_{\treen}[\mu_0; \boldsymbol{\varphi}]$ \emph{is a random pm, and}
\begin{equation} \label{eq:QCtnbis}
\mathcal{Q}_{\treen}[\mu_0](D) = \int_{[0, \pi]^{n-1}} \mathcal{C}_{\treen}[\mu_0; \boldsymbol{\varphi}](D) \beta^{\otimes_{n-1}}(\ud \boldsymbol{\varphi})
\end{equation}
\emph{holds true for every} $n \geq 2$ \emph{and} $D$ \emph{in} $\borelthree$.
\end{prop}

\emph{Proof}: As to the measurability of $\mathcal{C}_{\treen}[\mu_0; \boldsymbol{\varphi}]$, note that $\mathcal{C}_{\treen}$ gives a pm for every $\boldsymbol{\varphi}$ in $[0, \pi]^{n-1}$. Then, proceed to show that its Fourier transform is measurable as a function of $(\xib, \boldsymbol{\varphi})$. By resorting, once again, to Proposition 9 in Section 9.3 of \cite{frgr}, it suffices to verify that $\hat{\mathcal{C}}_{\treen}[\mu_0; \boldsymbol{\varphi}](\xib)$ is continuous as a function of $\xib$, for every fixed $\boldsymbol{\varphi}$, and measurable as a function of $\boldsymbol{\varphi}$, for every fixed $\xib$. The former claim is immediate in view of the basic properties of the Fourier transform, while the latter can be proved by induction. Fix $\xib \neq \mathbf{0}$ and combine (\ref{eq:bobylevforC}) with (\ref{eq:Cmisureiterate}) to write
\begin{eqnarray}
\hat{\mathcal{C}}_{\treen}[\mu_0; \boldsymbol{\varphi}](\xib) &=& \int_{0}^{2\pi} \hat{\mathcal{C}}_{\mathfrak{t}_{n}^{l}}[\mu_0; \boldsymbol{\varphi}^l](\rho \cos\varphi_{n-1} \boldsymbol{\psi}^l(\varphi_{n-1}, \theta, \ub))  \nonumber \\
&\times& \hat{\mathcal{C}}_{\mathfrak{t}_{n}^{r}}[\mu_0; \boldsymbol{\varphi}^r](\rho \sin\varphi_{n-1} \boldsymbol{\psi}^r( \varphi_{n-1}, \theta, \ub)) u_{(0, 2\pi)}(\ud \theta) \ . \ \ \ \ \ \ \ \ \ \label{eq:fourierCiterate}
\end{eqnarray}
Assume, as inductive hypothesis, the measurability of $\boldsymbol{\varphi} \mapsto \mathcal{C}_{\tree}[\mu_0; \boldsymbol{\varphi}]$ for every $\tree \in \cup_{\substack{h \leq n-1}} \mathbb{T}(h)$. Then, invoke the already mentioned proposition in \cite{frgr} to prove that both
$$
\begin{array}{l}
\hat{\mathcal{C}}_{\mathfrak{t}_{n}^{l}}[\mu_0; \boldsymbol{\varphi}^l](\rho \cos\varphi_{n-1} \boldsymbol{\psi}^l(\varphi_{n-1}, \theta, \ub)) \\
\hat{\mathcal{C}}_{\mathfrak{t}_{n}^{r}}[\mu_0; \boldsymbol{\varphi}^r](\rho \sin\varphi_{n-1} \boldsymbol{\psi}^r(\varphi_{n-1}, \theta, \ub))
\end{array}
$$
are measurable. Clearly, it is enough to study the former, the analysis of the latter being analogous. First, as $\rho$, $\ub$ and $\theta$ are fixed, think of $\hat{\mathcal{C}}_{\treen^l}[\mu_0; \boldsymbol{\varphi}^l](\rho \cos\varphi_{n-1} \boldsymbol{\psi}^l(\varphi_{n-1}, \theta,  \ub))$ as a function of the pair $(\boldsymbol{\varphi}^l, \varphi_{n-1})$ and note that the dependence on $\varphi_{n-1}$ is continuous, for any given $\boldsymbol{\varphi}^l$. To show this, it suffices to check the position of  $\varphi_{n-1}$ in the definition of $\boldsymbol{\psi}^l$ (see (\ref{eq:psi})) and to observe that $\rho \cos\varphi_{n-1} \boldsymbol{\psi}^l$ is the argument of a characteristic function. Next, for fixed $\varphi_{n-1}$, $\boldsymbol{\varphi}^l \mapsto \hat{\mathcal{C}}_{\treen^l}[\mu_0; \boldsymbol{\varphi}^l](\rho \cos\varphi_{n-1} \boldsymbol{\psi}^l(\varphi_{n-1}, \theta, \ub))$ is measurable thanks to the inductive hypothesis. This proves the measurability of the Fourier transform of $\mathcal{C}_{\treen}[\mu_0; \boldsymbol{\varphi}]$. At this stage, the same argument as in the proof of Proposition \ref{prop:C}, based on the L\'{e}vy inversion formula, leads to the measurability of $\boldsymbol{\varphi} \mapsto \mathcal{C}_{\treen}[\mu_0; \boldsymbol{\varphi}]$.

Finally, as far as (\ref{eq:QCtnbis}) is concerned, one first proves the equality
\begin{equation} \label{eq:QCtnbisfourier}
\hat{\mathcal{Q}}_{\treen}[\mu_0](\xib) = \int_{[0, \pi]^{n-1}} \hat{\mathcal{C}}_{\treen}[\mu_0; \boldsymbol{\varphi}](\xib) \beta^{\otimes_{n-1}}(\ud \boldsymbol{\varphi})
\end{equation}
for $n = 2, 3, \dots$, by mathematical induction. Indeed, (\ref{eq:QCtnbisfourier}) coincides with (\ref{eq:parametrizedbobylev}) when $n = 2$. For $n \geq 3$, combine (\ref{eq:Qtn}) with (\ref{eq:parametrizedbobylev}) to obtain
\begin{eqnarray}
\hat{\mathcal{Q}}_{\treen}[\mu_0](\xib) &=& \hat{\mathcal{Q}}\left[\mathcal{Q}_{\treen^l}[\mu_0], \mathcal{Q}_{\treen^r}[\mu_0]\right](\xib) \nonumber \\
&=& \int_{0}^{\pi}\int_{0}^{2\pi} \hat{\mathcal{Q}}_{\treen^l}[\mu_0](\rho \cos\varphi \boldsymbol{\psi}^l) \hat{\mathcal{Q}}_{\treen^r}[\mu_0](\rho \sin\varphi \boldsymbol{\psi}^r) u_{(0, 2\pi)}(\ud \theta) \beta(\ud \varphi) \ . \nonumber
\end{eqnarray}
By the inductive hypothesis,
\begin{eqnarray}
&{}& \int_{0}^{\pi}\int_{0}^{2\pi} \hat{\mathcal{Q}}_{\treen^l}[\mu_0](\rho \cos\varphi \boldsymbol{\psi}^l) \hat{\mathcal{Q}}_{\treen^r}[\mu_0](\rho \sin\varphi \boldsymbol{\psi}^r) u_{(0, 2\pi)}(\ud \theta) \beta(\ud \varphi) \nonumber \\
&=& \int_{0}^{2\pi}\int_{0}^{\pi} \int_{[0, \pi]^{n_l-1}} \hat{\mathcal{C}}_{\treen^l}[\mu_0; \boldsymbol{\varphi}^l](\rho \cos\varphi \boldsymbol{\psi}^l) \beta^{\otimes_{n_l-1}}(\ud \boldsymbol{\varphi}^l) \times \nonumber \\
&\times& \int_{[0, \pi]^{n_r-1}} \hat{\mathcal{C}}_{\treen^r}[\mu_0; \boldsymbol{\varphi}^r](\rho \sin\varphi \boldsymbol{\psi}^r) \beta^{\otimes_{n_r-1}}(\ud \boldsymbol{\varphi}^r) \beta(\ud \varphi) u_{(0, 2\pi)}(\ud \theta) \nonumber \\
&=& \int_{[0, \pi]^{n-1}}\int_{0}^{2\pi} \hat{\mathcal{C}}_{\treen^l}[\mu_0; \boldsymbol{\varphi}^l](\rho \cos\varphi \boldsymbol{\psi}^l) \hat{\mathcal{C}}_{\treen^r}[\mu_0; \boldsymbol{\varphi}^r](\rho \sin\varphi \boldsymbol{\psi}^r) \times \nonumber \\
&\times& u_{(0, 2\pi)}(\ud \theta) \beta^{\otimes_{n-1}}(\ud \boldsymbol{\varphi})
\ . \nonumber
\end{eqnarray}
This yields (\ref{eq:QCtnbisfourier}) by means of (\ref{eq:bobylevforC}) and (\ref{eq:Cmisureiterate}). Then, equality (\ref{eq:QCtnbis}) follows from (\ref{eq:QCtnbisfourier}) through Lemma \ref{lm:mixturechacteristics}.

\subsection{Description of the probabilistic framework} \label{sect:McKean}

The sample space $\Omega$, mentioned firstly in Theorem \ref{thm:main} and used throughout the work, is defined to be
$$
\Omega := \mathbb{N} \times \mathbb{T} \times [0, \pi]^{\infty} \times (0, 2\pi)^{\infty} \times (\rthree)^{\infty}
$$
where, for every nonempty set $X$, $X^{\infty}$ stands for the set of all sequences $(x_1, x_2, \dots)$ whose elements belong to $X$. The $\sigma$-algebra $\mathscr{F}$ is given by
$$
\mathscr{F} := \mathscr{B}(\mathbb{N}) \otimes  \mathscr{B}(\mathbb{T}) \otimes  \mathscr{B}([0, \pi]^{\infty}) \otimes  \mathscr{B}((0, 2\pi)^{\infty}) \otimes  \mathscr{B}((\rthree)^{\infty})
$$
where $\mathbb{N}$ and $\mathbb{T}$ are endowed with the discrete topology and $\mathscr{B}(X)$ stands for the Borel class in $X$. The symbols
$$
\nu,\ (\tau_n)_{n \geq 1},\ (\phi_n)_{n \geq 1},\ (\vartheta_n)_{n \geq 1}, (\mathbf{X}_n)_{n \geq 1}
$$
denote the coordinate random variables of $\Omega$. They can be given the following meaning: $\nu$ is the number of leaves, $(\tau_n)_{n \geq 1}$ is a sequence of trees with $\tau_n$ in $\mathbb{T}(n)$ for every $n \geq 1$, $(\phi_n)_{n \geq 1}$ and $(\vartheta_n)_{n \geq 1}$ are sequences of angles taking values in $[0, \pi]$ and $(0, 2\pi)$, respectively, and $(\mathbf{X}_n)_{n \geq 1}$ is a sequence of velocities. For each $t \geq 0$, $\pt$ is defined to be the pm on $(\Omega, \mathscr{F})$ which makes \emph{these random variables all stochastically independent of one another}, with these properties:
\begin{enumerate}
\item[a)] $\nu$ takes values in $\mathbb{N}$ and
\begin{equation} \label{eq:nu}
\pt[ \nu = n ] = e^{-t}(1 - e^{-t})^{n-1} \ \ \ \ \ \ \ (n = 1, 2, \dots)
\end{equation}
with the proviso that $0^0 := 1$.
\item[b)] $(\tau_n)_{n \geq 1}$ is a Markov sequence driven by
\begin{eqnarray}
\pt[\tau_1 = \tree_1] = 1 \label{eq:markov0} \\
\pt[ \tau_{n+1} = \treenk \ | \ \tau_n = \treen ] = \frac{1}{n} \label{eq:markov}
\end{eqnarray}
for every $n$, $\treen$ in $\mathbb{T}(n)$ and $k = 1, \dots, n$, and
$$
\pt[ \tau_{n+1} = \trees \ | \ \tau_n = \treen ] = 0
$$
whenever $\trees \not\in \mathbb{G}(\treen)$.
\item[c)] The elements of $(\phi_n)_{n \geq 1}$ are iid random numbers with common distribution $\beta$, specified in (\ref{eq:beta}).
\item[d)] The elements of $(\vartheta_n)_{n \geq 1}$ are iid with common distribution $u_{(0, 2\pi)}$.
\item[e)] As stated in Section \ref{sect:intro}, the $\mathbf{X}_j$s are iid with
\begin{equation} \label{eq:Xj}
\pt[\mathbf{X}_1 \in D] = \mu_0(D)
\end{equation}
for every $D$ in $\borelthree$.
\end{enumerate}
As to the existence of $\pt$ satisfying the above properties, see Theorem 3.19 in \cite{ka}.

There is an interesting relationship between the weights $p_n$, given in (\ref{eq:pntn}), and the law of the Markov sequence in b).
\begin{lm} \label{lm:pntn}
\emph{Equality}
\begin{equation}
p_n(\treen) = \pt[\tau_n = \treen]
\end{equation}
\emph{holds true for every} $n$ \emph{and} $\treen$ \emph{in} $\mathbb{T}(n)$.
\end{lm}

\emph{Proof}: Argue by mathematical induction. First, the assertion is trivially true for $n = 1, 2$. Next, given any $n \geq 3$,
\begin{eqnarray}
\pt[\tau_n = \treen] &=& \sum_{\treesn \in \ptreen} \pt[\tau_n = \treen \ | \ \tau_{n-1} = \treesn] \ \pt[\tau_{n-1} = \treesn] \nonumber \\
&=& \frac{1}{n-1} \sum_{\treesn \in \ptreen} \pt[\tau_{n-1} = \treesn] = \frac{1}{n-1} \sum_{\treesn \in \ptreen} p_{n-1}(\treesn) \nonumber
\end{eqnarray}
the last equality being valid thanks to the inductive hypothesis. Now,
$$
\sum_{\treesn \in \ptreen} p_{n-1}(\treesn) = \sum_{\substack{\treesn \in \ptreen \\
\treesn^l = \treen^l}} p_{n-1}(\treesn) \ + \sum_{\substack{\treesn \in \ptreen \\
\treesn^r = \treen^r}} p_{n-1}(\treesn)
$$
and, by (\ref{eq:pntn}),
\begin{eqnarray}
&{}& \sum_{\substack{\treesn \in \ptreen \\
\treesn^l = \treen^l}} p_{n-1}(\treesn) \ + \sum_{\substack{\treesn \in \ptreen \\
\treesn^r = \treen^r}} p_{n-1}(\treesn) \nonumber \\
&=& \frac{1}{n-2} \Big{[} \sum_{\substack{\treesn \in \ptreen \\
\treesn^l = \treen^l}} p_{n_l}(\treesn^l) p_{n_r - 1}(\treesn^r) \ + \sum_{\substack{\treesn \in \ptreen \\
\treesn^r = \treen^r}} p_{n_l - 1}(\treesn^l) p_{n_r}(\treesn^r) \Big{]} \nonumber \\
&=& \frac{1}{n-2} \Big{[} p_{n_l}(\treen^l) \sum_{\substack{\treesn \in \ptreen \\
\treesn^l = \treen^l}} p_{n_r - 1}(\treesn^r) \ + \ p_{n_r}(\treen^r) \sum_{\substack{\treesn \in \ptreen \\
\treesn^r = \treen^r}} p_{n_l - 1}(\treesn^l) \Big{]} \ . \nonumber
\end{eqnarray}
At this stage, one can conclude by working out the previous sums as follows. As to the sum relative to $\treesn^l = \treen^l$,
\begin{eqnarray}
&{}& \sum_{\substack{\treesn \in \ptreen \\
\treesn^l = \treen^l}}
p_{n_r - 1}(\treesn^r) = \sum_{\mathfrak{s}_{n_r - 1} \in \mathbb{P}(\treen^r)} p_{n_r - 1}(\mathfrak{s}_{n_r - 1}) \nonumber \\
&=& \sum_{\mathfrak{s}_{n_r - 1} \in \mathbb{P}(\treen^r)} \pt[\tau_{n_r - 1} = \mathfrak{s}_{n_r - 1}] \nonumber \\
&=& (n_r - 1) \sum_{\mathfrak{s}_{n_r - 1} \in \mathbb{P}(\treen^r)} \pt[\tau_{n_r - 1} = \mathfrak{s}_{n_r - 1}] \pt[\tau_{n_r} = \treen^r \ | \ \tau_{n_r - 1} = \mathfrak{s}_{n_r - 1}] \nonumber \\
&=& (n_r - 1) \pt[\tau_{n_r} = \treen^r] = (n_r - 1) p_{n_r}(\treen^r) \nonumber \ .
\end{eqnarray}
The same line of reasoning leads to
$$
\sum_{\substack{\treesn \in \ptreen \\ \treesn^r = \treen^r}} p_{n_l - 1}(\treesn^l) = (n_l - 1) p_{n_l}(\treen^l) \nonumber \ .
$$
The proof of the lemma terminates by combining the previous equations, upon noting that $n = n_l + n_r$. \ $\Box$ \\

The way is paved for the definition of the arrays $(\pi_{1, n}, \dots, \pi_{n, n})$ and $(\mathrm{O}_{1, n}, \dots, \mathrm{O}_{n, n})$ mentioned in Section \ref{sect:intro}. As to the former, let $\pi_{j, n}^{\ast}$ be a real-valued function on $\mathbb{T}(n) \times [0, \pi]^{n-1}$, for $j = 1, \dots, n$ and $n$ in $\mathbb{N}$. Specifically, $\pi_{1, 1}^{\ast} \equiv 1$ and, for $n \geq 2$,
\begin{equation} \label{eq:pijnast}
\pi_{j, n}^{\ast}(\treen, \boldsymbol{\varphi}) := \left\{ \begin{array}{ll}
\pi_{j, n_l}^{\ast}(\mathfrak{t}_{n}^{l}, \boldsymbol{\varphi}^l) \cos\varphi_{n-1} & \text{for} \ j = 1, \dots, n_l \\
\pi_{j - n_l, n - n_l}^{\ast}(\mathfrak{t}_{n}^{r}, \boldsymbol{\varphi}^r) \sin\varphi_{n-1} & \text{for} \ j = n_l + 1, \dots, n
\end{array} \right.
\end{equation}
for every $\boldsymbol{\varphi} = (\boldsymbol{\varphi}^l, \boldsymbol{\varphi}^r, \varphi_{n-1})$ in $[0, \pi]^{n-1}$. At this stage, set
\begin{equation} \label{eq:pijn}
\pi_{j, n} := \pi_{j, n}^{\ast}(\tau_n, (\phi_1, \dots, \phi_{n-1}))
\end{equation}
for $j = 1, \dots, n$ and $n$ in $\mathbb{N}$. A straightforward induction argument shows that
\begin{equation} \label{eq:sumpijn}
\sum_{j = 1}^{n} \pi_{j, n}^2 = 1
\end{equation}
for every $n$ in $\mathbb{N}$ and $\omega$ in $\Omega$. By arguing as in \cite{blm, gr6}, verify that
\begin{equation} \label{eq:gare}
\et\left[\sum_{j = 1}^{\nu} |\pi_{j, \nu}|^s \right] = e^{-(1 - 2 l_s(b))t}
\end{equation}
holds true for every $s > 0$ and $l_s(b) := \int_{0}^{\pi} (\sin\varphi)^s \beta(\ud \varphi)$.

The definition of the $\mathrm{O}_{j, n}$s requires the introduction of the following $3 \times 3$ orthogonal matrices
$$
\mathrm{M}^l(\varphi, \theta) := \left( \begin{array}{ccc} -\cos\theta \cos\varphi & \sin\theta & \cos\theta \sin\varphi \\
- \sin\theta \cos\varphi & -\cos\theta & \sin\theta \sin\varphi  \\
\sin\varphi & 0 & \cos\varphi \\
\end{array} \right)
$$
and
$$
\mathrm{M}^r(\varphi, \theta) := \left( \begin{array}{ccc} \sin\theta & \cos\theta \sin\varphi & -\cos\theta \cos\varphi \\
-\cos\theta & \sin\theta \sin\varphi & - \sin\theta \cos\varphi \\
0 & \cos\varphi & \sin\varphi \\
\end{array} \right)
$$
on the basis of which one considers $\mathbb{SO}(3)$-valued functions $\mathrm{O}_{j, n}^{\ast}$ on $\mathbb{T}(n) \times [0, \pi]^{n-1} \times (0, 2\pi)^{n-1}$, for $j = 1, \dots, n$ and $n$ in $\mathbb{N}$. They are defined by
$$
\mathrm{O}_{1, 1}^{\ast} \equiv \mathrm{Id}_{3 \times 3}
$$
and, for $n \geq 2$,
\begin{eqnarray} \label{eq:Ojn}
&{}& \mathrm{O}_{j, n}^{\ast}(\treen, \boldsymbol{\varphi}, \boldsymbol{\theta}) \nonumber \\
&:=& \left\{ \begin{array}{ll}
\mathrm{M}^l(\varphi_{n-1}, \theta_{n-1}) \mathrm{O}_{j, n_l}^{\ast}(\treen^l, \boldsymbol{\varphi}^l, \boldsymbol{\theta}^l) & \text{for} \ j = 1, \dots, n_l \\
\mathrm{M}^r(\varphi_{n-1}, \theta_{n-1}) \mathrm{O}_{j - n_l, n - n_l}^{\ast}(\treen^r, \boldsymbol{\varphi}^r, \boldsymbol{\theta}^r) & \text{for} \ j = n_l + 1, \dots, n
\end{array} \right.
\end{eqnarray}
for every $\boldsymbol{\varphi}$ in $[0, \pi]^{n-1}$ and $\boldsymbol{\theta}$ in $(0, 2\pi)^{n-1}$. Development of this recursion formula gives
\begin{equation} \label{eq:pathOijn}
\mathrm{O}_{j, n}^{\ast}(\treen, \boldsymbol{\varphi}, \boldsymbol{\theta}) = \prod_{h = 1}^{\delta_j(\treen)} \mathrm{M}^{\epsilon_h(\treen, j)}(\varphi_{m_h(\treen, j)}, \theta_{m_h(\treen, j)})
\end{equation}
where $\prod_{h=1}^{n} A_h := A_1 \times \dots \times A_n$. The $\epsilon_h(\treen, j)$s take values in $\{l, r\}$ and, in particular, $\epsilon_h(\treen, j)$ equals $l$ ($r$, respectively) if $j \leq n_l$ ($j > n_l$, respectively). Each $m_h$ belongs to $\{1, \dots, n-1\}$ and $m_1 \neq \dots \neq m_{\delta_j(\treen)}$. In particular, $m_1(\treen, j) = n-1$, for every $\treen$ in $\mathbb{T}(n)$ and $j = 1, \dots, n$. It is worth mentioning that there is a more direct procedure to prove      (\ref{eq:pathOijn}) based on the specific structure of each $\treen$. It is explained in detail in Section 2.5 of \cite{dophd}.

The random matrices mentioned in Section \ref{sect:intro} can now be specified by
\begin{equation} \label{eq:Ojn}
\mathrm{O}_{j, n} := \mathrm{O}_{j, n}^{\ast}(\tau_n, (\phi_1, \dots, \phi_{n-1}), (\vartheta_1, \dots, \vartheta_{n-1}))
\end{equation}
for $j = 1, \dots, n$ and $n$ in $\mathbb{N}$. To conclude, note that $\psib_{j, n}$, defined by (\ref{eq:psijn}), is a random function of the variable $\ub$ satisfying
\begin{equation} \label{eq:moduluspsi}
|\psib_{j, n}(\ub)| = 1
\end{equation}
for $\ub$ in $S^2$, $j = 1, \dots, n$ and $n$ in $\mathbb{N}$.

\subsection{Proof of Theorem \ref{thm:main}} \label{sect:proof1}

Substitute (\ref{eq:QCtnbis}) in (\ref{eq:McKean}) to have
\begin{eqnarray}
\mu(D, t) &=& e^{-t}\mu_0(D) + \sum_{n = 2}^{+\infty} e^{-t} (1 - e^{-t})^{n-1} \sum_{\treen \in \mathbb{T}(n)} p_n(\treen) \times \nonumber \\
&\times& \int_{[0, \pi]^{n-1}} \mathcal{C}_{\treen}[\mu_0; \boldsymbol{\varphi}](D) \beta^{\otimes_{n-1}}(\ud \boldsymbol{\varphi}) \label{eq:CDGRmeas}
\end{eqnarray}
for every $D$ in $\borelthree$ and every $t \geq 0$. Next, define the function $\mathcal{M} : \Omega \rightarrow \prob$ as
\begin{equation} \label{eq:randmeas}
\mathcal{M}(D) := \mathcal{C}_{\tau_{\nu}}[\mu_0; (\phi_1, \dots, \phi_{\nu - 1})](D)
\end{equation}
for every $D$ in $\borelthree$, with the proviso that $(\phi_1, \dots, \phi_{\nu - 1}) := \emptyset$ if $\nu = 1$.
In view of Proposition \ref{prop:Ctn}, $\mathcal{M}$ is $\mathscr{F}/\mathscr{B}(\prob)$-measurable, which is tantamount to saying that $\mathcal{M}$ is a random pm. Moreover, $\mathcal{M}$ does not depend on $b$. Taking expectation of both sides of (\ref{eq:randmeas}) yields
\begin{eqnarray}
\et[\mathcal{M}(D)] &=& \int_{\Omega} \mathcal{C}_{\tau_{\nu}}[\mu_0; (\phi_1, \dots, \phi_{\nu - 1})](D) \ud\pt \nonumber \\
&=& e^{-t}\mu_0(D) + \sum_{n = 2}^{+\infty} \sum_{\treen \in \mathbb{T}(n)} \ \int_{[0, \pi]^{n-1}} \mathcal{C}_{\treen}[\mu_0; (\varphi_1, \dots, \varphi_{n-1})](D) \times \nonumber \\
&\times& \pt[\nu = n, \tau_n = \treen, \phi_1 \in \ud \varphi_1, \dots, \phi_{n - 1} \in \ud \varphi_{n-1}] \ . \label{eq:expectationM}
\end{eqnarray}
After recalling the assumption of independence of the coordinate random variables, combination of points a)-c) in Subsection \ref{sect:McKean} with Lemma \ref{lm:pntn} gives
\begin{gather}
\pt[\nu = n, \tau_n = \treen, \phi_1 \in \ud \varphi_1, \dots, \phi_{n - 1} \in \ud \varphi_{n-1}] \nonumber \\
= e^{-t} (1 - e^{-t})^{n-1}  p_n(\treen)  \beta^{\otimes_{n-1}}(\ud \varphi_1, \dots, \ud \varphi_{n-1}) \nonumber
\end{gather}
for every $n \geq 2$, which shows that the RHSs in (\ref{eq:CDGRmeas}) and (\ref{eq:expectationM}) coincide.

\subsection{Proof of Theorem \ref{thm:Mhat}} \label{sect:proof2}

After defining the $\sigma$-algebras
\begin{eqnarray}
\mathscr{G} &:=& \sigma\big{(}\nu,\ (\tau_n)_{n \geq 1},\ (\phi_n)_{n \geq 1} \big{)} \nonumber \\
\mathscr{H} &:=& \sigma\big{(}\nu,\ (\tau_n)_{n \geq 1},\ (\phi_n)_{n \geq 1},\ (\vartheta_n)_{n \geq 1} \big{)} \nonumber
\end{eqnarray}
it is easy to verify points from i) to v) and the validity of the equalities
\begin{eqnarray}
\hat{\mathcal{M}}(\xib) &=& \hat{\mathcal{C}}_{\tau_{\nu}}[\mu_0; (\phi_1, \dots, \phi_{\nu - 1})] \nonumber \\
\textsf{E}_t \left[e^{i \rho S(\ub)} \ | \ \mathscr{H} \right] &=& \textsf{E}_t \left[\prod_{j = 1}^{\nu} \hat{\mu}_0(\rho \pi_{j, \nu} \psib_{j, \nu}) \ \big{|} \ \mathscr{H} \right] \nonumber \\
\textsf{E}_t \left[e^{i \rho S(\ub)} \ | \ \mathscr{G} \right] &=& \textsf{E}_t \left[\prod_{j = 1}^{\nu} \hat{\mu}_0(\rho \pi_{j, \nu} \psib_{j, \nu}) \ \big{|} \ \mathscr{G} \right] \nonumber
\end{eqnarray}
so that the problem reduces to proving that
\begin{eqnarray}
&{}& \hat{\mathcal{C}}_{\treen}[\mu_0; \boldsymbol{\varphi}](\xib) \nonumber \\
&=& \int_{(0, 2\pi)^{n-1}} \left[\prod_{j=1}^{n} \hat{\mu}_0\big{(} \rho  \pi_{j, n}^{\ast}(\treen, \boldsymbol{\varphi})
\qb_{j, n}(\treen, \boldsymbol{\varphi}, \boldsymbol{\theta}, \ub)
\big{)}\right] u_{(0, 2\pi)}^{\otimes_{n-1}}(\ud \boldsymbol{\theta}) \ \ \ \ \ \ \label{eq:FourierCtnast} \\
&=& \int_{(0, 2\pi)^{n-1}} \left[\prod_{j=1}^{n} \hat{\mu}_0\big{(}\rho \pi_{j, n}^{\ast}(\treen, \boldsymbol{\varphi})
\mathrm{B}(\ub) \mathrm{O}_{j, n}^{\ast}(\treen, \boldsymbol{\varphi}, \boldsymbol{\theta}) \mathbf{e}_3
\big{)}\right] u_{(0, 2\pi)}^{\otimes_{n-1}}(\ud \boldsymbol{\theta}) \ \ \ \ \ \ \label{eq:FourierCtn}
\end{eqnarray}
hold true for every $n \geq 2$, for every tree $\treen$ in $\mathbb{T}(n)$, $\boldsymbol{\varphi}$ in $[0, \pi]^{n-1}$ and $\xib \neq \mathbf{0}$, with $\rho = |\xib|$ and $\ub = \xib/|\xib|$. In (\ref{eq:FourierCtnast}), the $\qb_{j, n}$s are defined as follows. First, $\qb_{1, 1}(\tree_1, \emptyset, \emptyset, \ub) := \ub$. Then, for every tree $\treen$ with $n \geq 2$, put
\begin{eqnarray} \label{eq:psijnabstract}
&{}& \qb_{j, n}(\treen, \boldsymbol{\varphi}, \boldsymbol{\theta}, \ub) \nonumber \\
&=& \left\{ \begin{array}{ll}
\qb_{j, n_l}(\treen^l, \boldsymbol{\varphi}^l, \boldsymbol{\theta}^l, \psib^l(\varphi_{n-1}, \theta_{n-1}, \ub)) & \text{for} \ j = 1, \dots, n_l \\
\qb_{j - n_l, n - n_l}(\treen^r, \boldsymbol{\varphi}^r, \boldsymbol{\theta}^r, \psib^r(\varphi_{n-1}, \theta_{n-1}, \ub)) & \text{for} \ j = n_l + 1, \dots, n
\end{array} \right.
\end{eqnarray}
for every $\boldsymbol{\varphi} = (\boldsymbol{\varphi}^l, \boldsymbol{\varphi}^r, \varphi_{n-1})$ in $[0, \pi]^{n-1}$ and $\boldsymbol{\theta} = (\boldsymbol{\theta}^l, \boldsymbol{\theta}^r, \theta_{n-1})$ in $(0, 2\pi)^{n-1}$, where $\psib^l$ and $\psib^r$ are given by (\ref{eq:psi}).

To prove (\ref{eq:FourierCtnast}), first consider the case when $n = 2$ and observe that $\pi_{1, 2}^{\ast} = \cos\varphi_1$, $\pi_{2, 2}^{\ast} = \sin\varphi_1$, $\qb_{1, 2} = \psib^l$, $\qb_{2, 2} = \psib^r$. Then, check that (\ref{eq:FourierCtnast}) reduces to (\ref{eq:bobylevforC}) with $\zeta = \eta = \mu_0$. Next, by mathematical induction, assume $n \geq 3$ and recall (\ref{eq:bobylevforC}) again and (\ref{eq:Cmisureiterate}) to write
\begin{eqnarray}
\hat{\mathcal{C}}_{\treen}[\mu_0; \boldsymbol{\varphi}](\boldsymbol{\xi}) &=& \int_{0}^{2\pi} \hat{\mathcal{C}}_{\treen^l}[\mu_0; \boldsymbol{\varphi}^l](\rho \cos\varphi_{n-1} \boldsymbol{\psi}^l(\varphi_{n-1}, \theta_{n-1}, \ub))  \nonumber \\
&\times& \hat{\mathcal{C}}_{\treen^r}[\mu_0; \boldsymbol{\varphi}^r](\rho \sin\varphi_{n-1} \boldsymbol{\psi}^r(\varphi_{n-1}, \theta_{n-1}, \ub)) u_{(0, 2\pi)}(\ud \theta_{n-1}) \ . \nonumber
\end{eqnarray}
Thus, assuming that (\ref{eq:FourierCtnast}) holds true for every $l$ in $\{1, \dots, n-1\}$ and every tree $\tree_{l}$ in $\cup_{l \leq n-1} \mathbb{T}(l)$, deduce that
\begin{gather}
\hat{\mathcal{C}}_{\treen^l}[\mu_0; \boldsymbol{\varphi}^l](\rho \cos\varphi_{n-1} \psib^l(\varphi_{n-1}, \theta_{n-1}, \ub)) \nonumber \\
\begin{array}{lll}
&= \int_{(0, 2\pi)^{n_l - 1}} \Big{\{}&\prod_{j=1}^{n_l} \hat{\mu}_0 \big{[} \rho \cos\varphi_{n-1} \pi_{j, n_l}^{\ast}(\treen^l, \boldsymbol{\varphi}^l) \\
&{} &\times \qb_{j, n_l}(\treen^l, \boldsymbol{\varphi}^l, \boldsymbol{\theta}^l, \boldsymbol{\psi}^l(\varphi_{n-1}, \theta_{n-1}, \ub))\big{]} \Big{\}} u_{(0, 2\pi)}^{\otimes_{n_l - 1}}(\ud \boldsymbol{\theta}^l) \nonumber
\end{array}
\end{gather}
and
\begin{gather}
\hat{\mathcal{C}}_{\treen^r}[\mu_0; \boldsymbol{\varphi}^r](\rho \cos\varphi_{n-1} \boldsymbol{\psi}^r(\varphi_{n-1}, \theta_{n-1}, \ub)) \nonumber \\
\begin{array}{lll}
&= \int_{(0, 2\pi)^{n_r - 1}} \Big{\{}&\prod_{j=1}^{n_r} \hat{\mu}_0 \big{[} \rho \cos\varphi_{n-1} \pi_{j, n_r}^{\ast}(\treen^r, \boldsymbol{\varphi}^r) \\
&{} &\times \qb_{j, n_r}(\treen^r, \boldsymbol{\varphi}^r, \boldsymbol{\theta}^r, \boldsymbol{\psi}^r(\varphi_{n-1}, \theta_{n-1}, \ub))\big{]} \Big{\}} u_{(0, 2\pi)}^{\otimes_{n_r - 1}}(\ud \boldsymbol{\theta}^r) \ . \nonumber
\end{array}
\end{gather}
To complete the argument, combine the above equalities, invoke (\ref{eq:pijnast}) and (\ref{eq:psijnabstract}) and note that $u_{(0, 2\pi)} \otimes u_{(0, 2\pi)}^{\otimes_{n_l - 1}} \otimes u_{(0, 2\pi)}^{\otimes_{n - n_l - 1}} = u_{(0, 2\pi)}^{\otimes_{n - 1}}$.

As far as the proof of (\ref{eq:FourierCtn}) is concerned, start by noting that
$$
\qb_{j, 2}(\tree_2, \varphi, \theta, \ub) = \mathrm{B}(\ub)\mathrm{O}_{j, 2}^{\ast}(\tree_2, \varphi, \theta) \mathbf{e}_3
$$
for $j = 1, 2$, for every $\varphi$ in $[0, \pi]$, $\theta$ in $(0, 2\pi)$ and $\ub$ in $S^2$, provided that the basis $\{\mathbf{a}(\ub), \mathbf{b}(\ub), \ub\}$ in (\ref{eq:psi}) is formed by the three columns of $\mathrm{B}(\ub)$.
Then, assume $n \geq 3$ and argue by induction, starting from (\ref{eq:FourierCtnast}) and definitions (\ref{eq:Ojn}) and (\ref{eq:psijnabstract}). Whence,
\begin{gather}
\int_{(0, 2\pi)^{n-1}} \left[\prod_{j=1}^{n} \hat{\mu}_0\big{(} \rho  \pi_{j, n}^{\ast}(\treen, \boldsymbol{\varphi}) \qb_{j, n}(\treen, \boldsymbol{\varphi}, \boldsymbol{\theta}, \ub)\big{)}\right] u_{(0, 2\pi)}^{\otimes_{n-1}}(\ud \boldsymbol{\theta}) \nonumber \\
= \int_{0}^{2\pi} \! \int_{(0, 2\pi)^{n_l-1}} \int_{(0, 2\pi)^{n_r-1}} \nonumber \\
\left[\prod_{j=1}^{n_l} \hat{\mu}_0\big{(} \rho  \pi_{j, n}^{\ast}(\treen, \boldsymbol{\varphi}) \qb_{j, n_l}(\treen^l, \boldsymbol{\varphi}^l, \boldsymbol{\theta}^l, \boldsymbol{\psi}^l(\varphi_{n-1}, \theta_{n-1}, \ub))\big{)}\right] \times \nonumber \\
\times  \left[\prod_{j=1}^{n_r} \hat{\mu}_0\big{(} \rho  \pi_{j + n_l, n}^{\ast}(\treen, \boldsymbol{\varphi}) \qb_{j, n_r}(\treen^r, \boldsymbol{\varphi}^r, \boldsymbol{\theta}^r, \boldsymbol{\psi}^r(\varphi_{n-1}, \theta_{n-1}, \ub))\big{)}\right] \times \nonumber \\
\times u_{(0, 2\pi)}^{\otimes_{n_r - 1}}(\ud \boldsymbol{\theta}^r)  u_{(0, 2\pi)}^{\otimes_{n_l - 1}}(\ud \boldsymbol{\theta}^l) u_{(0, 2\pi)}(\ud  \theta_{n-1}) \nonumber \\
= \int_{0}^{2\pi} \! \int_{(0, 2\pi)^{n_l-1}} \int_{(0, 2\pi)^{n_r-1}} \nonumber \\
\left[\prod_{j=1}^{n_l} \hat{\mu}_0\big{(} \rho  \pi_{j, n}^{\ast}(\treen, \boldsymbol{\varphi}) \mathrm{B}(\psib^l(\varphi_{n-1}, \theta_{n-1}, \ub)) \mathrm{O}_{j, n_l}^{\ast}(\treen^l, \boldsymbol{\varphi}^l, \boldsymbol{\theta}^l) \mathbf{e}_3 \big{)}\right] \times \nonumber \\
\times  \left[\prod_{j=1}^{n_r} \hat{\mu}_0\big{(} \rho  \pi_{j + n_l, n}^{\ast}(\treen, \boldsymbol{\varphi})
\mathrm{B}(\psib^r(\varphi_{n-1}, \theta_{n-1}, \ub)) \mathrm{O}_{j, n_r}^{\ast}(\treen^r, \boldsymbol{\varphi}^r, \boldsymbol{\theta}^r) \mathbf{e}_3 \big{)}\right] \times \nonumber \\
\times u_{(0, 2\pi)}^{\otimes_{n_r - 1}}(\ud \boldsymbol{\theta}^r)  u_{(0, 2\pi)}^{\otimes_{n_l - 1}}(\ud \boldsymbol{\theta}^l) u_{(0, 2\pi)}(\ud  \theta_{n-1})  \ . \label{eq:intermediatepsiastpsi}
\end{gather}
The integral $\int_{(0, 2\pi)^{n_l-1}}$ ($\int_{(0, 2\pi)^{n_l-1}}$, respectively) should not be written if $n_l = 1$ ($n_r = 1$, respectively) since $\boldsymbol{\theta}^l$ ($\boldsymbol{\theta}^r$, respectively) corresponds to the empty set. At this stage, it will be proved that
\begin{gather}
\int_{(0, 2\pi)^{n_l-1}} \left[\prod_{j=1}^{n_l} \hat{\mu}_0\big{(} x \mathrm{B}(\psib^l(\varphi_{n-1}, \theta_{n-1}, \ub)) \mathrm{O}_{j, n_l}^{\ast}(\treen^l, \boldsymbol{\varphi}^l, \boldsymbol{\theta}^l) \mathbf{e}_3 \big{)}\right]  u_{(0, 2\pi)}^{\otimes_{n_l - 1}}(\ud \boldsymbol{\theta}^l) \nonumber \\
= \int_{(0, 2\pi)^{n_l-1}} \left[\prod_{j=1}^{n_l} \hat{\mu}_0\big{(} x
\mathrm{B}(\ub) \mathrm{O}_{j, n}^{\ast}(\treen, \boldsymbol{\varphi}, \boldsymbol{\theta}) \mathbf{e}_3
\big{)}\right] u_{(0, 2\pi)}^{\otimes_{n_l - 1}}(\ud \boldsymbol{\theta}^l) \label{eq:robben}
\end{gather}
for every real $x$, $\ub$ in $S^2$, $\boldsymbol{\varphi}$ in $[0, \pi]^{n-1}$ and $\theta_{n-1}$ in $(0, 2\pi)$. If $n_l = 1$, it is immediate to verify that
\begin{equation} \label{eq:orthogonaltrick}
^t\psib^l(\varphi_{n-1}, \theta_{n-1}, \ub) = \mathrm{B}(\psib^l(\varphi_{n-1}, \theta_{n-1}, \ub)) \mathbf{e}_3 = \mathrm{B}(\ub) \mathrm{M}^l(\varphi_{n-1}, \theta_{n-1}) \mathbf{e}_3
\end{equation}
which implies (\ref{eq:robben}). Then, pass to the most interesting case when $n_l \geq 2$. In view of (\ref{eq:orthogonaltrick}), the last column of the two matrices $\mathrm{B}(\boldsymbol{\psi}^l(\varphi_{n-1}, \theta_{n-1}, \ub))$ and $\mathrm{B}(\ub) \mathrm{M}^l(\varphi_{n-1}, \theta_{n-1})$ are the same and, therefore, there exists an orthogonal matrix of the type
$$
\mathrm{R}(\alpha) := \left( \begin{array}{ccc}
\cos\alpha & -\sin\alpha & 0 \\
\sin\alpha & \cos\alpha & 0 \\
0 & 0 & 1 \\
\end{array} \right)
$$
for which $\mathrm{B}(\boldsymbol{\psi}^l(\varphi_{n-1}, \theta_{n-1}, \ub)) = \mathrm{B}(\ub) \mathrm{M}^l(\varphi_{n-1}, \theta_{n-1}) \mathrm{R}(\alpha)$. Note also that $\alpha$ depends only on $\varphi_{n-1}$, $\theta_{n-1}$ and $\ub$. Now, when $n_l \geq 2$, recall that
$$
\mathrm{O}_{j, n_l}(\treen^l, \boldsymbol{\varphi}^l, \boldsymbol{\theta}^l) = \mathrm{M}^{\epsilon_{1}(\treen^l, j)}(\varphi_{n_l-1}, \theta_{n_l-1}) \mathrm{Q}_{j, n_l}
$$
where $\mathrm{Q}_{j, n_l}$ is another orthogonal matrix. If $\delta_j(\treen^l) = 1$, $\mathrm{Q}_{j, n_l}$ reduces to the identity matrix. Otherwise, it depends on $\treen^l$, $(\varphi_1, \dots, \varphi_{n_l -2})$ and
$(\theta_1, \dots, \theta_{n_l -2})$ through (\ref{eq:pathOijn}). In any case, the precise expression is not needed here. The key remark is represented by the equality
$$
\mathrm{R}(\alpha) \mathrm{M}^s(\varphi, \theta) = \mathrm{M}^s(\varphi, \theta + \alpha)
$$
which holds for both $s = l, r$ and for every $\varphi$ and $\theta$. Then, start from the LHS of (\ref{eq:robben}) and take account of the product $\mathrm{R}(\alpha)\mathrm{M}^{\epsilon_{1}(\treen^l, j)}(\varphi_{n_l-1}, \theta_{n_l-1})$, which equals $\mathrm{M}^{\epsilon_{1}(\treen^l, j)}(\varphi_{n_l-1}, \theta_{n_l-1} + \alpha(\varphi_{n-1}, \theta_{n-1}, \ub))$. A change of variable $\theta_{n_l-1}^{'} = \theta_{n_l-1} + \alpha(\varphi_{n-1}, \theta_{n-1}, \ub)$ transforms the LHS of (\ref{eq:robben}) into
$$
\int_{(0, 2\pi)^{n_l-1}} \left[\prod_{j=1}^{n_l} \hat{\mu}_0\big{(} x \mathrm{B}(\ub) \mathrm{M}^l(\varphi_{n-1}, \theta_{n-1}) \mathrm{O}_{j, n_l}^{\ast}(\treen^l, \boldsymbol{\varphi}^l, \boldsymbol{\theta}^l) \mathbf{e}_3 \big{)}\right]  u_{(0, 2\pi)}^{\otimes_{n_l - 1}}(\ud \boldsymbol{\theta}^l)
$$
which equals the RHS of (\ref{eq:robben}), in view of (\ref{eq:Ojn}).

To complete the proof of (\ref{eq:FourierCtn}), use (\ref{eq:intermediatepsiastpsi}) after noting that an equality similar to (\ref{eq:robben}) can be stated by changing subscripts and superscripts from $l$ to $r$, and replacing $\mathrm{O}_{j + n_l, n}^{\ast}$ to $\mathrm{O}_{j, n}^{\ast}$.

Finally, the above argument works in the same way to show that two representations of the type of (\ref{eq:FourierCtn}) are equivalent even if definition (\ref{eq:psijn}) is based, through the matrix $\mathrm{B}(\ub)$, on two different choices of orthonormal bases, such as $\{\mathbf{a}(\ub), \mathbf{b}(\ub), \ub\}$ and $\{\mathbf{a}^{'}(\ub), \mathbf{b}^{'}(\ub), \ub\}$.

\section{Central limit theorem for Maxwellian molecules} \label{sect:CLT}

We are going to give a necessary and sufficient condition for the convergence of the solution $\mu(\cdot, t)$ of (\ref{eq:boltzmann}), as $t$ goes to infinity, by resorting to arguments which are reminiscent of those used to prove central limit theorems. This is made possible in view of Theorems \ref{thm:main}-\ref{thm:Mhat}, which express $\mu$ in terms of sums of random numbers. As recalled in Section \ref{sect:intro}, the study of the convergence to equilibrium of the solution of (\ref{eq:boltzmann}) has been considered as a problem of paramount interest ever since the pioneering works in kinetic theory of gases. As a consequence, this question has given rise to an extensive research. See Chapter 2D in \cite{vil}.

Before stating the main result of this section, it is worth recalling that (\ref{eq:secondmoment}) implies
\begin{eqnarray}
\intethree \vb \mu(\ud \vb, t) &=& \intethree \vb \mu_0(\ud \vb) \label{eq:Vu0} \\
\intethree |\vb|^2 \mu(\ud \vb, t) &=& \intethree |\vb|^2 \mu_0(\ud \vb) \label{eq:sigma0}
\end{eqnarray}
for every $t \geq 0$. Moreover, the unique stationary solutions of (\ref{eq:boltzmann}) -- in the class of all pms on $\borelthree$ -- are the so-called \emph{Maxwellian distributions} given by
\begin{equation} \label{eq:maxwellians}
\gamma_{\vb_0, \sigma}(\ud \vb) = \left(\frac{1}{2\pi \sigma^2}\right)^{3/2} \exp\{-\frac{1}{2 \sigma^2} |\vb - \vb_0|^2\} \ud \vb
\end{equation}
with $(\vb_0, \sigma)$ in $\rthree \times [0, +\infty)$. If $\sigma = 0$, $\gamma_{\vb_0, 0}$ reduces to the unit mass $\delta_{\vb_0}$ at $\vb_0$. We are now in a position to formulate the foretold result.
\begin{thm} \label{thm:CLT}
\emph{Let conditions} (\ref{eq:maxwell})-(\ref{eq:cutoff}) \emph{be in force. Then}, $\mu(\cdot, t)$ \emph{converges weakly, as} $t$ \emph{goes to infinity, if and only if} (\ref{eq:secondmoment}) \emph{holds true. In case this condition is satisfied, the limiting distribution is given by} (\ref{eq:maxwellians}) \emph{with}
\begin{eqnarray}
\vb_0 &=& \intethree \vb \mu_0(\ud \vb) \nonumber \\
\sigma^2 &=& \frac{1}{3} \intethree |\vb - \vb_0|^2 \mu_0(\ud \vb) \ . \nonumber
\end{eqnarray}
\end{thm}

At the beginning, starting from initial data with finite second moment and satisfying extra conditions related to the Boltzmann H-theorem, different authors stated the convergence to the Maxwellian distribution with respect to the total variation distance. Afterwards, on the basis of inspiring principles different from the H-theorem, the aforesaid conclusions
have been improved in different ways. In \cite{ta} it has been proved that the solution \emph{converges weakly} when (\ref{eq:secondmoment}) is valid. A definitive result in \cite{cl} has finally stated that (\ref{eq:secondmoment}) suffices to obtain convergence in total variation. To the best of the authors' knowledge, the only proof that (\ref{eq:secondmoment}) is necessary is contained in \cite{cgr}, where it has been shown that if the second moment of $\mu_0$ is infinite then all the mass ``explodes to infinity''. The motivations for the new proof, we are going to present, rely on the nexus with the methods of current usage for solving the central limit problem, which one can now institute thanks to Theorems \ref{thm:main}-\ref{thm:Mhat}. Indeed, the cohesive power of these methods gives rise to a well-structured, transparent and direct proof that turns out to be, on the whole, more comprehensible than a mere union of necessarily disparate arguments extracted from the above-mentioned works. An analogous strategy has already been followed in \cite{gr8} apropos of the simpler case of Kac's equation.

\subsection{Sufficiency of condition (\ref{eq:secondmoment})} \label{sect:sufficient}

If $\intethree |\xb|^2 \mu_0(\ud \xb) = 0$, then $\mu(\cdot, t)$ coincides with the unit mass at $\mathbf{0}$, $\delta_{\mathbf{0}}(\cdot)$, for every $t \geq 0$. On the other hand, when (\ref{eq:secondmoment}) is in force with $\intethree |\xb|^2 \mu_0(\ud \xb) > 0$, one can assume that the following extra-conditions are valid, without any loss of generality.
\begin{enumerate}
\item[i)] $\intethree \xb \mu_0(\ud \xb) = \mathbf{0}$
\item[ii)] $\intethree |\xb|^2 \mu_0(\ud \xb) = 3$
\item[iii)] $\intethree x_i x_j \mu_0(\ud \xb) = 0$, whenever $1 \leq i \neq j \leq 3$
\item[iv)] $\sigma_{\ast}^2 := \min\{\sigma_{1}^2, \sigma_{2}^2, \sigma_{3}^2\} > 0$, with $\sigma_{i}^2 := \intethree x_{i}^{2} \mu_0(\ud \xb)$ for $i = 1, 2, 3$.
\end{enumerate}
In fact, i) and ii) can be assumed in view of the conservation property encapsulated in (\ref{eq:Vu0})-(\ref{eq:sigma0}). Under i) and ii), the thesis of Theorem \ref{thm:CLT} can be specified by saying that the parameters of the limiting Maxwellian are $\vb_0 = \mathbf{0}$ and $\sigma^2 = 1$. Condition iii) can be assumed in view of Lemma \ref{lm:orthogonal}, by choosing $R$ in such a way that $\ ^tR V[\mu_0] R$ is a diagonal matrix, where $V[\mu_0]$ denotes the covariance matrix of $\mu_0$. At this stage, replacement of $\mu_0$ with $\mu_0 \circ \mathcal{R}^{-1}$, where $\mathcal{R} : \xb \mapsto R\xb$, does not alter the conclusion to be reached. Note that the ensuing initial datum fulfills conditions i)-iii). As to condition iv), recall that $\intethree x_{i}^{2} \mu(\ud \xb, t) \rightarrow 1$ as $t$ goes to infinity, for $i = 1, 2, 3$. See Lemma 1.8 in \cite{ccg0}. Hence, for each strictly positive $\eta$, there is $t_0 = t_0(\eta, b)$ so that
\begin{equation} \label{eq:variance1}
\big{|} \intethree x_{i}^{2} \mu(\ud \xb, t) - 1 \big{|} \leq \eta
\end{equation}
for every $t \geq t_0$, and $i = 1, 2, 3$. Substituting $\mu_0$ with $\mu(\cdot, t_0)$, iv) holds true along with i)-iii).

The sufficiency of (\ref{eq:secondmoment}) will be now proved by means of the L\'{e}vy continuity theorem for characteristic functions, by adapting the argument displayed, for example, in Section 9.1 of \cite{chte}. To start, fix $\xib \neq \mathbf{0}$, put $\rho := |\xib|$ and $\ub := \xib/\rho$. Theorem \ref{thm:main} gives
\begin{eqnarray}
\big{|}\hat{\mu}(\xib, t) - \exp\{-|\xib|^2/2\}\big{|} &=& \big{|}\et\hat{\mathcal{M}}(\xib) - \exp\{-|\xib|^2/2\}\big{|} \nonumber \\
&\leq& \big{|}\et[\hat{\mathcal{M}}(\rho \ub) - \exp\{-T^2 \rho^2/2\}]\big{|} \nonumber \\
&+& \big{|}\et\exp\{-T^2 \rho^2/2\} - \exp\{-\rho^2/2\}\big{|} \label{eq:suff1}
\end{eqnarray}
where
$$
T^2 := \et\left[\Big{(}S(\ub)\Big{)}^2\ \Big{|} \ \mathscr{H}\right] = \sum_{j =1}^{\nu} \pi_{j, \nu}^2 \sum_{i = 1}^3 \sigma_{i}^2 \psi_{j, \nu; i}^2(\ub) \ .
$$
As to the first summand in (\ref{eq:suff1}), Theorem \ref{thm:Mhat} can be invoked to write
\begin{equation} \label{eq:suff6}
\big{|}\et[\hat{\mathcal{M}}(\rho \ub) - \exp\{-T^2 \rho^2/2\}]\big{|} \leq \et \Big{|} \et\left[e^{i \rho S(\ub)}\ | \ \mathscr{H}\right] - \exp\{-T^2 \rho^2/2\}\Big{|} \ .
\end{equation}
By using (6)-(7) in Section 9.1 of \cite{chte}, one gets
\begin{eqnarray}
&{}& \Big{|}\et\left[e^{i \rho S(\ub)}\ | \ \mathscr{H}\right] - \exp\{-T^2 \rho^2/2\}\Big{|} \nonumber \\
&\leq& \rho^2 \sum_{j = 1}^{\nu} \pi_{j, \nu}^2 \int_{A_j(\varepsilon)} (\psib_{j, \nu}(\ub) \cdot \xb)^2 \mu_0(\ud \xb)
+ \varepsilon |T|^3 \rho^3 \nonumber \\
&+& \frac{1}{8} T^4 \rho^4 \frac{\max_{1 \leq j \leq \nu} \pi_{j, \nu}^2 \sum_{i = 1}^3 \sigma_{i}^2 \psi_{j, \nu; i}^2(\ub)}{T^2} \label{eq:suff5}
\end{eqnarray}
with
$$
A_j(\varepsilon) := \{\xb \in \rthree \ \big{|} \ |\pi_{j, \nu} (\psib_{j, \nu}(\ub) \cdot \xb)| \geq \varepsilon |T|\}
$$
for $j = 1, \dots, \nu$. Combination of points i)-iv) with (\ref{eq:sumpijn}) and (\ref{eq:moduluspsi}) shows that $\sigma_{\ast}^2 \leq T^2 \leq 3$ and
\begin{equation} \label{eq:suff3}
\frac{1}{8} T^2 \rho^4 \max_{1 \leq j \leq \nu} \pi_{j, \nu}^2 \sum_{i = 1}^3 \sigma_{i}^2 \psi_{j, \nu; i}^2(\ub) \leq \frac{9}{8} \rho^4 \pi_{o}^2
\end{equation}
with $\pi_o := \max_{1 \leq j \leq \nu} |\pi_{j, \nu}|$. Put $M(y) := \int_{|\xb| \geq 1/y} |\xb|^2 \mu_0(\ud \xb)$ for every $y$ in $(0, +\infty)$ and note that $M$ is a monotonically increasing, bounded function satisfying $\lim_{y \downarrow 0} M(y) = 0$. Moreover, from
\begin{eqnarray}
A_j(\varepsilon) &\subset& \{\xb \in \rthree\ \big{|} \ |\pi_{j, \nu} (\psib_{j, \nu}(\ub) \cdot \xb)| \geq \varepsilon \sigma_{\ast}\} \nonumber \\
&\subset& \{\xb \in \rthree\ \big{|} \ |\pi_o (\psib_{j, \nu}(\ub) \cdot \xb)| \geq \varepsilon \sigma_{\ast}\}
\subset \{\xb \in \rthree\ \big{|} \ |\pi_o| \cdot |\xb| \geq \varepsilon \sigma_{\ast}\} \nonumber
\end{eqnarray}
one can conclude that
\begin{equation} \label{eq:suff4}
\sum_{j = 1}^{\nu} \pi_{j, \nu}^2 \int_{A_j(\varepsilon)} (\psib_{j, \nu}(\ub) \cdot \xb)^2 \mu_0(\ud \xb) \leq M\left(\frac{\pi_o}{\varepsilon \sigma_{\ast}}\right)
\end{equation}
holds true for every strictly positive $\varepsilon$. At this stage, fix $\varepsilon = \sqrt{\pi_o}$ and combine (\ref{eq:suff5}) with (\ref{eq:suff3}) and (\ref{eq:suff4}) to get
\begin{equation} \label{eq:suff2}
\Big{|}\et\left[e^{i \rho S(\ub)}\ | \ \mathscr{H}\right] - \exp\{-T^2 \rho^2/2\}\Big{|} \leq M\left(\frac{\sqrt{\pi_o}}{\sigma_{\ast}}\right)\rho^2  + 3\sqrt{3} \sqrt{\pi_o} \rho^3 + \frac{9}{8} \pi_{o}^2 \rho^4 \ .
\end{equation}
To complete the analysis of the RHS of (\ref{eq:suff6}), it can be shown that the expectation of the RHS of (\ref{eq:suff2}) goes to zero as $t$ goes to infinity, for every $\rho$ in $(0, +\infty)$. Indeed, for any monotonically increasing, bounded function $g : (0, \infty) \rightarrow (0, \infty)$ satisfying $\lim_{x \downarrow 0} g(x) = 0$, one has
\begin{eqnarray}
\et[g(\pi_{o})] &=& \et[g(\pi_{o})\ind\{\pi_{o} \leq e^{-z t}\}] + \et[g(\pi_{o})\ind\{\pi_{o} > e^{-z t}\}] \nonumber \\
&\leq& g(e^{-z t}) + \sup_{x \in (0, \infty)} g(x) \cdot \pt[\pi_{o} > e^{-z t}] \nonumber \\
&\leq& g(e^{-z t}) + \sup_{x \in (0, \infty)} g(x) \cdot \et[\pi_{o}^4] e^{4 z t} \nonumber
\end{eqnarray}
for any $z$ in $(0, \infty)$. Using (\ref{eq:gare}) with $s = 4$ gives $\et[\pi_{o}^4] \leq e^{-(1 - 2l_4(\beta))t}$ and, after choosing $z = \frac{1}{8}(1 - 2l_4(b))$, one obtains $\lim_{t \rightarrow +\infty} \et[g(\pi_{o})] = 0$. This argument, with $g(x) = M\left(\frac{\sqrt{x}}{\sigma_{\ast}}\right)\rho^2  + 3\sqrt{3} \sqrt{x} \rho^3 + \frac{9}{8} x^2 \rho^4$, leads to the desired result.

As far as the second summand in (\ref{eq:suff1}) is concerned, one first observes that, in view of Theorem \ref{thm:main}, $\et[e^{- T^2 \rho^2/2}]$ can be thought of as the Fourier transform of the solution of (\ref{eq:boltzmann}) when the initial datum $\mu_0$ meets
\begin{equation} \label{eq:ingauss}
\mu_0(\ud \xb) = \prod_{i = 1}^{3} \frac{1}{\sigma_i \sqrt{2\pi}} \exp\{-\frac{x_{i}^{2}}{2 \sigma_{i}^{2}}\} \ud x_i \ .
\end{equation}
By resorting to an easy argument developed in \cite{doLom}, $\eta$ in $(\ref{eq:variance1})$ can be chosen so small that (\ref{eq:ingauss}) belongs to a convenient neighborhood of $\gamma_{\mathbf{0}, 1}$, so that
$$
\big{|}\et\exp\{-T^2 \rho^2/2\} - \exp\{-\rho^2/2\}\big{|} \leq c_1 e^{-c_2 t}
$$
for some strictly positive constants $c_1$, $c_2$. See Theorem 1 in \cite{doLom}.

\subsection{Necessity of condition (\ref{eq:secondmoment})} \label{sect:necessity}

Suppose that the solution $\mu(\cdot, t)$ of (\ref{eq:boltzmann}), with initial datum $\mu_0$, converges weakly to some limit, as $t$ goes to infinity. Following a technique developed in \cite{flr}, the argument starts with the introduction of the random vector
$$
W = \big{(}\nu, (\tau_j)_{j \geq 1}, (\phi_j)_{j \geq 1}, (\vartheta_j)_{j \geq 1}, \boldsymbol{\lambda}, \Lambda, \mathbf{U}\big{)}
$$
defined on $(\Omega, \mathscr{F})$. The meaning of $\nu, (\tau_j)_{j \geq 1}, (\phi_j)_{j \geq 1}, (\vartheta_j)_{j \geq 1}$ is the same as in Subsection \ref{sect:McKean} while, to explain the other symbols, one fixes an arbitrary point $\ub_0$ in $S^2$ and defines:
\begin{enumerate}
\item[i)] $\boldsymbol{\lambda} := (\lambda_1(\cdot), \dots, \lambda_{\nu}(\cdot), \delta_0(\cdot), \delta_0(\cdot), \dots)$ to be the sequence of random pms with Fourier transforms $\hat{\lambda}_j(\xi) := \hat{\mu}_0(\xi \pi_{j, \nu} \psib_{j, \nu}(\ub_0))$, for $j = 1, \dots, \nu$ and $\xi$ in $\rone$, while $\delta_x$ stands for the degenerate pm at $x$.
\item[ii)] $\Lambda$ to be the random pm obtained as convolution of all the elements of $\boldsymbol{\lambda}$, i.e. $\Lambda = \lambda_1 \ast \dots \ast \lambda_{\nu}$.
\item[iii)] $\mathbf{U} := (U_1, U_2, \dots)$ to be the sequence of random numbers defined by $U_k := \max_{1 \leq j \leq \nu} \lambda_j\left(\left[-\frac{1}{k}, \frac{1}{k}\right]^c\right)$ for every $k$ in $\mathbb{N}$.
\end{enumerate}
To grasp the usefulness of $W$, one can note that its components are the essential ingredients of the central limit problem for \emph{independent uniformly asymptotically negligible} summands. See Sections 16.6-9 of \cite{frgr}. Apropos of the negligibility condition, it is easy to prove that
$$
\lim_{t \rightarrow +\infty} \pt[U_k > \alpha] = 0
$$
holds for every $k$ in $\mathbb{N}$ and for every $\alpha$ in $(0, +\infty)$. In fact, observe that $|\psib_{j, \nu}| = 1$ entails
$\{\xb \in \rthree \ \big{|} \ |\pi_{j, \nu} \psib_{j, \nu} \cdot \xb| \geq 1/k\} \subset \{\xb \in \rthree \ \big{|} \ |\pi_{j, \nu} \xb| \geq 1/k\}$ and hence
$$
\{U_k > \alpha\} \subset \left\{\left[\max_{1 \leq j \leq \nu} \mu_0\{\xb \in \rthree \ \big{|} \ |\pi_{j, \nu} \xb| \geq 1/k\}\right] \geq \alpha\right\} \ .
$$
Then, apply the argument used to prove Lemma 2 in \cite{gr8}.

Now, think of the range of $W$ as a subset of
$$
\mathbb{S} := \overline{\mathbb{N}} \times \overline{\mathbb{T}} \times [0, \pi]^{\infty} \times [0, 2\pi]^{\infty} \times (\mathscr{P}(\overline{\mathbb{R}}))^{\infty} \times \mathscr{P}(\overline{\mathbb{R}}) \times [0, 1]^{\infty}
$$
where $\overline{\mathbb{T}}$ is the one-point compactification of the discrete topological space $\mathbb{T}$; $\overline{\mathbb{N}} := \{1, 2, \dots, +\infty\}$, $\overline{\mathbb{R}} := [-\infty, +\infty]$. According to results explained in Section 6.II of \cite{parth}, $\mathscr{P}(\overline{\mathbb{R}})$ can be metrized, consistently with the topology of weak convergence, in such a way that it turns out to be a separable, compact and complete metric space. It follows that $\mathbb{S}$ is a separable, compact and complete metric space with respect to the product topology and so the family of probability distributions $(\pt \circ W^{-1})_{t \geq 0}$ is tight. This implies that any sequence $(\textsf{P}_{t_m} \circ W^{-1})_{m \geq 1}$, such that $t_m$ strictly increases to infinity as $m$ goes to infinity, contains a subsequence $(\textsf{Q}_n)_{n \geq 1}$, with $\textsf{Q}_n := \textsf{P}_{t_{m_n}} \circ W^{-1}$, which converges weakly to a pm $\textsf{Q}$. It is worth noting that, thanks to the weak convergence of $\mu(\cdot, t)$, $\textsf{Q}$ is supported by
$$
\{+\infty\} \times \overline{\mathbb{T}} \times [0, \pi]^{\infty} \times [0, 2\pi]^{\infty} \times (\delta_0, \delta_0, \dots) \times \mathscr{P}(\mathbb{R}) \times \{0\}^{\infty} \ .
$$
This claim can be verified by invoking Lemma 3 in \cite{gr8}.

Since $\mathbb{S}$ is Polish, from the \emph{Skorokhod representation theorem} (see Theorem 6.7 in \cite{bill2}) one can determine a probability space $(\tilde{\Omega}, \tilde{\mathscr{F}}, \tilde{\pp})$ and random elements on it, taking values
in $\mathbb{S}$,
\begin{eqnarray}
\tilde{W}_n &=& \big{(}\tilde{\nu}_n, (\tilde{\tau}_{j, n})_{j \geq 1}, (\tilde{\phi}_{j, n})_{j \geq 1}, (\tilde{\vartheta}_{j, n})_{j \geq 1}, \tilde{\boldsymbol{\lambda}}_n, \tilde{\Lambda}_n, \tilde{\mathbf{U}}_n\big{)} \nonumber \\
\tilde{W}_{\infty} &=& \big{(}\{+\infty\}, (\tilde{\tau}_{j, \infty})_{j \geq 1}, (\tilde{\phi}_{j, \infty})_{j \geq 1}, (\tilde{\vartheta}_{j, \infty})_{j \geq 1}, (\delta_0, \delta_0, \dots), \tilde{\Lambda}_{\infty}, (0, 0, \dots)\big{)} \nonumber
\end{eqnarray}
which have respective probability laws $\textsf{Q}_n$ and $\textsf{Q}$ and satisfy $\tilde{W}_n(\tilde{\omega}) \rightarrow \tilde{W}_{\infty}(\tilde{\omega})$ in the metric of $\mathbb{S}$, as $n \rightarrow +\infty$, for every $\tilde{\omega}$ in $\tilde{\Omega}$. This entails
\begin{equation} \label{eq:convergenzeflr}
\begin{array}{lll}
\tilde{\nu}_n \rightarrow +\infty, & \ \ \ \tilde{\mathbf{U}}_n \rightarrow (0, 0, \dots) \\
\tilde{\boldsymbol{\lambda}}_n \Rightarrow (\delta_0, \delta_0, \dots), & \ \ \ \tilde{\Lambda}_n \Rightarrow \tilde{\Lambda}_{\infty}
\end{array}
\end{equation}
as $n$ goes to infinity. The symbol $\Rightarrow$ is here used to designate weak convergence of pms. The distributional properties of $\tilde{W}_n$ imply that $\tilde{\Lambda}_n$ is the convolution of the elements of $\tilde{\boldsymbol{\lambda}}_n$, and that $\tilde{U}_{k, n} = \max_{1 \leq j \leq \tilde{\nu}_n} \tilde{\lambda}_{j, n}\left(\left[-\frac{1}{k}, \frac{1}{k}\right]^c\right)$ for every $k$ in $\mathbb{N}$, $\tilde{\pp}$-almost surely.

From now on, given any pm $q$ on $\mathscr{B}(\rone^d)$, $q^{(s)}$ will denote the symmetrized $q$, i.e. $\hat{q}^{(s)}(\cdot) := |\hat{q}(\cdot)|^2$. Then, (\ref{eq:convergenzeflr}) entails $\tilde{\Lambda}_{n}^{(s)} \Rightarrow \tilde{\Lambda}_{\infty}^{(s)}$ for every $\tilde{\omega}$ in $\tilde{\Omega}$, which, combined with Theorem 24 in Chapter 16 of \cite{frgr}, gives
\begin{equation} \label{eq:frgr}
+\infty > \sigma^2(\tilde{\omega}) := \lim_{\varepsilon \downarrow 0} \overline{\underline{\lim}}_n \sum_{j = 1}^{\tilde{\nu}_n(\tilde{\omega})} \int_{[-\varepsilon, \varepsilon]} x^2 \tilde{\lambda}_{j}^{(s)}(\ud x; \tilde{\omega})
\end{equation}
with the exception of a set of points $\tilde{\omega}$ of $\tilde{\pp}$-probability 0. The proof of the necessity starts from this inequality and carries on according to the following argument. Firstly,
\begin{eqnarray}
\sum_{j = 1}^{\tilde{\nu}_n} \int_{[-\varepsilon, \varepsilon]} x^2 \tilde{\lambda}_{j}^{(s)}(\ud x) &=& \sum_{j = 1}^{\tilde{\nu}_n} \tilde{\pi}_{j, \tilde{\nu}_n}^{2} \intethree (\tilde{\psib}_{j, \tilde{\nu}_n} \cdot \xb)^2 \ind\{|\tilde{\pi}_{j, \tilde{\nu}_n} \tilde{\psib}_{j, \tilde{\nu}_n} \cdot \xb| \leq \varepsilon\} \mu_{0}^{(s)}(\ud \xb) \nonumber \\
&\geq& \sum_{j = 1}^{\tilde{\nu}_n} \tilde{\pi}_{j, \tilde{\nu}_n}^{2} \sum_{i = 1}^{3} \tilde{\psi}_{j, \tilde{\nu}_n; i}^{2} \int_{\tilde{\pi}_{n}^{\ast} |\xb| \leq \varepsilon} x_{i}^{2} \mu_{0}^{(s)}(\ud \xb)
\end{eqnarray}
where $\tilde{\pi}$ and $\tilde{\psib}$ denote the counterparts, in the Skorokhod representation, of the $\pi$ and $\psib$ defined in Subsection \ref{sect:McKean}, $\tilde{\pi}_{n}^{\ast} := \max_{1 \leq j \leq \tilde{\nu}_n} |\tilde{\pi}_{j, \tilde{\nu}_n}|$ and the inequality is a consequence of the inclusion
$$
\{\xb \in \rthree \ \big{|} \ |\tilde{\pi}_{j, \tilde{\nu}_n} \tilde{\psib}_{j, \tilde{\nu}_n} \cdot \xb| \leq \varepsilon\} \supset \{\xb \in \rthree \ \big{|} \ \tilde{\pi}_{n}^{\ast} |\xb| \leq \varepsilon\} \ .
$$
Secondly, define $k = k(\tilde{\omega}; j, n)$ to be an element of $\{1, 2, 3\}$ for which $\tilde{\psi}_{j, \tilde{\nu}_n; k}^2 = \max_{1 \leq i \leq 3} \tilde{\psi}_{j, \tilde{\nu}_n; i}^2$, which is greater than $1/3$ since $\tilde{\psib}_{j, \tilde{\nu}_n}$ belongs to $S^2$, for every $\tilde{\omega}$ in $\tilde{\Omega}$, $n$ in $\mathbb{N}$ and $j = 1, \dots, \tilde{\nu}_n$. Then,
\begin{eqnarray}
\sum_{j = 1}^{\tilde{\nu}_n} \tilde{\pi}_{j, \tilde{\nu}_n}^{2} \sum_{i = 1}^{3} \tilde{\psi}_{j, \tilde{\nu}_n; i}^2 \int_{\tilde{\pi}_{n}^{\ast} |\xb| \leq \varepsilon} x_{i}^{2} \mu_{0}^{(s)}(\ud \xb) &\geq& \sum_{j = 1}^{\tilde{\nu}_n} \tilde{\pi}_{j, \tilde{\nu}_n}^{2} \tilde{\psi}_{j, \tilde{\nu}_n; k}^2 \int_{\tilde{\pi}_{n}^{\ast} |\xb| \leq \varepsilon} x_{k}^{2} \mu_{0}^{(s)}(\ud \xb) \nonumber \\
&\geq& \frac{1}{3} \sum_{j = 1}^{\tilde{\nu}_n} \tilde{\pi}_{j, \tilde{\nu}_n}^{2} \int_{\tilde{\pi}_{n}^{\ast} |\xb| \leq \varepsilon} x_{k}^{2} \mu_{0}^{(s)}(\ud \xb) \nonumber \\
&=& \frac{1}{3} \sum_{h = 1}^{3} \tilde{s}_{h, n} \int_{\tilde{\pi}_{n}^{\ast} |\xb| \leq \varepsilon} x_{h}^{2} \mu_{0}^{(s)}(\ud \xb) \nonumber
\end{eqnarray}
where $\tilde{s}_{h, n}$ denotes the sum of those $\tilde{\pi}_{j, \tilde{\nu}_n}^{2}$ so that $k(\tilde{\omega}; j, n) = h$. In view of Lemma 1 in \cite{gr8}, without real loss of generality one can assume that $\tilde{\pi}_{n}^{\ast}$ goes to zero with probability one, as $n$ goes to infinity. Hence, since $\sum_{h = 1}^{3} \tilde{s}_{h, n} = 1$ with probability one, there are some $\tilde{\omega}$ and $h$, say $h_{\ast} = h_{\ast}(\tilde{\omega})$, so that $\tilde{\pi}_{n}^{\ast}(\tilde{\omega}) \rightarrow 0$ and $\overline{\lim}_n \tilde{s}_{h_{\ast}, n}$ is strictly positive. Then,
\begin{eqnarray}
\sigma^2(\tilde{\omega}) &\geq& \lim_{\varepsilon \downarrow 0} \overline{\lim}_n \frac{1}{3} \sum_{h = 1}^{3} \tilde{s}_{h, n} \int_{\tilde{\pi}_{n}^{\ast} |\xb| \leq \varepsilon} x_{h}^{2} \mu_{0}^{(s)}(\ud \xb) \nonumber \\
&\geq& \frac{1}{3} \overline{\lim}_{r \rightarrow +\infty} \int_{|\xb| \leq r} x_{h_{\ast}}^{2} \mu_{0}^{(s)}(\ud \xb) \cdot \overline{\lim}_n \tilde{s}_{h_{\ast}, n} \nonumber
\end{eqnarray}
which is enough to show that the $h_{\ast}$-th marginal of $\mu_{0}^{(s)}$ -- and hence also the $h_{\ast}$-th marginal of $\mu_{0}$ -- has finite second moment.

To complete the proof, recall Lemma \ref{lm:orthogonal} and note that weak convergence of $\mu(\cdot, t)$ entails weak convergence of $\mu(\cdot, t) \circ \mathcal{R}^{-1}$. Hence, the above argument can be used to prove that $\intethree x_{h_{\ast}}^{2} \mu_0 \circ \mathcal{R}^{-1}(\ud \xb)$ is finite for every choice of the orthogonal matrix $R$, with the same $h_{\ast}$, independently of $R$ and $\mu_0$. At this stage, it suffices to choose $\mathcal{R}$ firstly equal to $\mathcal{R}_1 := (x_1, x_2, x_3) \mapsto (x_2, x_3, x_1)$ and, then, equal to $\mathcal{R}_2 := (x_1, x_2, x_3) \mapsto (x_3, x_1, x_2)$ to have the desired result.

\section{Global boundedness of the moments of $\mathcal{M}$} \label{sect:elmroth}

The present section deals with the moments 
$$
\mathfrak{M}_r := \intethree |\xb|^{r} \mathcal{M}(\ud \xb)
$$
of the random pm $\mathcal{M}$. Since $\mathcal{M}$ depends on $\mu_0$, it is interesting to investigate the connection between the finiteness of $\mathfrak{M}_r$ and the finiteness of the moments of $\mu_0$. This problem is intimately linked with the global boundedness of the moments of the solution of (\ref{eq:boltzmann}), through Theorems \ref{thm:main}-\ref{thm:Mhat}. Global boundedness of the moments of $\mu$ has been carefully investigated in some papers. See \cite{elm} and the references therein.

As an application of Theorems \ref{thm:main}-\ref{thm:Mhat}, we are now able to determine upper bounds, which are expressed as functions of the moments of $\mu_0$ and are independent of $b$, for the random moments $\mathfrak{M}_r$. To appreciate the novelty of this statement, recall that existing literature is limited to inequalities concerning only the expectation of the above random moments.

\begin{prop}
\emph{Assume that} $\intethree \xb \mu_0(\ud \xb) = \mathbf{0}$, $\intethree |\xb|^2 \mu_0(\ud \xb) = 3$ \emph{and that}
$$
\mathfrak{m}_r := \intethree |\xb|^r \mu_0(\ud \xb) < +\infty
$$
\emph{for} $r = 3, \dots, 2k$ \emph{and some integer} $k \geq 2$. \emph{Then, the random number} $S(\ub)$ \emph{defined by} (\ref{eq:Su}) \emph{satisfies}
\begin{equation} \label{eq:momlS}
\et\left[|S(\ub)|^l \ | \ \mathscr{H} \right] \leq g_l \\
\end{equation}
\emph{for} $l = 2, \dots, 2k$ \emph{and for every} $\ub$ \emph{in} $S^2$, $\pt$-\emph{almost surely. The upper bounds} $g_l$ \emph{are positive constants depending on} $\mu_0$ \emph{only through} $\mathfrak{m}_l$. \emph{Consequently}
\begin{equation} \label{eq:boundderivativeMhat}
\Big{|}\frac{\partial^l}{\partial \rho^l} \hat{\mathcal{M}}(\rho \ub)\Big{|} \leq g_l
\end{equation}
\emph{for} $l = 2, \dots, 2k$, $\rho$ \emph{in} $\rone$ \emph{and} $\ub$ \emph{in} $S^2$, $\pt$-\emph{almost surely. Moreover,}
\begin{equation} \label{eq:fullmoment}
\mathfrak{M}_{2h} \leq 3^h g_{2h}
\end{equation}
\emph{for} $h = 1, \dots, k$, $\pt$-\emph{almost surely.}
\end{prop}
The expression of the constants $g_l$ can be found in the following

\emph{Proof}: First, note that
$$
\et\left[(S(\ub))^2 \ | \ \mathscr{H} \right] = \sum_{j=1}^{\nu} \pi_{j, \nu}^2 \et\left[(\mathbf{X}_j \cdot \boldsymbol{\psi}_{j, \nu})^2 \ | \ \mathscr{H} \right] \leq 3
$$
as an obvious consequence of the conditional Cauchy-Schwarz inequality. Whence, (\ref{eq:momlS}) holds true for $l = 2$ with $g_2 := 3$. Then, for $l \geq 3$, an inequality due to Rosenthal -- see Section 2.3 in \cite{pe1} -- yields
\begin{gather}
\et\left[|S(\ub)|^l \ | \ \mathscr{H} \right] \nonumber \\
\leq c(l) \Big{\{} \sum_{j = 1}^{\nu} \et\left[ |\pi_{j, \nu} \mathbf{X}_j \cdot \boldsymbol{\psi}_{j, \nu}|^l \ | \ \mathscr{H} \right] + \Big{(} \sum_{j = 1}^{\nu} \et\left[|\pi_{j, \nu} \mathbf{X}_j \cdot \boldsymbol{\psi}_{j, \nu}|^2\ | \ \mathscr{H} \right] \Big{)}^{l/2} \Big{\}}  \nonumber
\end{gather}
where $c(l)$ is a positive constant depending only on $l$. An additional application of the Cauchy-Schwarz inequality, combined with (\ref{eq:sumpijn}) and (\ref{eq:moduluspsi}), gives
\begin{eqnarray}
\et\left[|S(\ub)|^l \ | \ \mathscr{H} \right] &\leq& c(l) \Big{\{} \mathfrak{m}_l \sum_{j = 1}^{\nu} |\pi_{j, \nu}|^l + \Big{(} 3 \sum_{j = 1}^{\nu} \pi_{j, \nu}^2 \Big{)}^{l/2}\Big{\}} \nonumber \\
&=& c(l) \Big{\{} \mathfrak{m}_l \sum_{j = 1}^{\nu} |\pi_{j, \nu}|^l +  3^{l/2}\Big{\}} \nonumber
\end{eqnarray}
which entails (\ref{eq:momlS}) with $g_l := c(l) \{\mathfrak{m}_l  +  3^{l/2}\}$. By elementary properties of the characteristic function, (\ref{eq:momlS}) entails
$$
\Big{|}\frac{\partial^l}{\partial \rho^l} \et\left[e^{i \rho S(\ub)} \ | \ \mathscr{H} \right] \Big{|} \leq g_l
$$
for $l = 2, \dots, 2k$, which gives (\ref{eq:boundderivativeMhat}) through an application of the dominated convergence theorem in (\ref{eq:Mhat}). As for (\ref{eq:fullmoment}), verify that
\begin{eqnarray}
\mathfrak{M}_{2h} &=& \intethree \left(\sum_{i = 1}^{3} x_{i}^2\right)^h \mathcal{M}(\ud \xb) \leq 3^{h-1} \sum_{i = 1}^{3}
\intethree x_{i}^{2h} \mathcal{M}(\ud \xb) \nonumber \\
&=& 3^{h-1} \sum_{i = 1}^{3} \lim_{\rho \rightarrow 0} \frac{\partial^{2h}}{\partial \rho^{2h}} \hat{\mathcal{M}}(\rho \mathbf{e}_i) \ , \nonumber
\end{eqnarray}
where $\mathbf{e}_i$ stands for the $i$-th vector of the canonical basis of $\rthree$, and conclude by applying (\ref{eq:boundderivativeMhat}). \\

\textbf{Acknowledgments.} Emanuele Dolera and Eugenio Regazzini are grateful to Eric Carlen for several useful discussions. They are partially supported by Italian Ministry of University and Research, grant no. 2008MK3AFZ. Eugenio Regazzini is affiliated also with IMATI-CNR, Milano, Italy.

\vspace{1cm}

\footnotesize{\textsc{emanuele dolera \\
dipartimento di matematica pura e applicata ``giuseppe vitali'' \\
universit\`a degli studi di modena e reggio emilia \\
via campi 213/b, 41100 modena, italy \\
e-mail:} emanuele.dolera@unimore.it, emanuele.dolera@unipv.it}

\vspace{1cm}

\footnotesize{\textsc{eugenio regazzini \\
dipartimento di matematica ``felice casorati'' \\
universit\`a degli studi di pavia \\
via ferrata 1, 27100 pavia, italy \\
e-mail:} eugenio.regazzini@unipv.it}


\begin{thebibliography}{99}

\bibitem{blm} \textsc{Bassetti, F., Ladelli, L.} and \textsc{Matthes, D.} (2010). Central limit theorem for a class of one-dimensional kinetic equations. \emph{Probab. Theory Relat. Fields}. DOI 10.1007/s00440-010-0269-8.

\bibitem{blr} \textsc{Bassetti, F., Ladelli, L.} and \textsc{Regazzini, E.} (2008). Probabilistic study of the speed of approach to equilibrium for an inelastic Kac model. \emph{J. Stat. Phys}. $\mathbf{133}$ 683-710.

\bibitem{bill2} \textsc{Billingsley, P.} (1999). \emph{Convergence of Probability Measures}. $2^{nd}$ ed. Wiley, NY.

\bibitem{bob88} \textsc{Bobylev, A.V.} (1988). The theory of the nonlinear spatially uniform Boltzmann equation for Maxwell molecules. \emph{Mathematical Physics Reviews}. $\mathbf{7}$ 111-233.

\bibitem{cc} \textsc{Carlen, E. A.} and \textsc{Carvalho, M. C.} (2003). Probabilistic methods in kinetic theory. \emph{Riv. Mat. Univ. Parma}. $\textbf{7}$ ($2^{\ast}$)101-149.

\bibitem{ccg0} \textsc{Carlen, E. A., Carvalho, M. C.} and \textsc{Gabetta, E.} (2000). Central limit theorem for Maxwellian molecules and truncation of the Wild expansion. \emph{Comm. Pure Appl. Math}. \textbf{53} 370-397.

\bibitem{cgr} \textsc{Carlen, E.A., Gabetta, E.} and \textsc{Regazzini, E.} (2007). On the rate of explosion for infinite energy solutions of the spatially homogeneous Boltzmann equation. \emph{J.
Stat. Phys.} $\mathbf{129}$ 699-723.

\bibitem{cl} \textsc{Carlen, E.A.} and \textsc{Lu, X.} (2003). Fast and slow convergence to
equilibrium for Maxwellian molecules via Wild sums. \emph{J.
Stat. Phys.} $\mathbf{112}$ 59-134.

\bibitem{cerS} \textsc{Cercignani, C.} (1988). \emph{The Boltzmann Equation and its Applications}. Springer-Verlag, NY.

\bibitem{chte} \textsc{Chow, Y. S.} and \textsc{Teicher, H.} (1997). \emph{Probability Theory. Independence, Interchangeability, Martingales}. $3^{rd}$ ed. Springer-Verlag, NY.

\bibitem{desv} \textsc{Desvillettes, L.} (2003). About the use of the Fourier transform for the Boltzmann equation. \emph{Riv. Mat. Univ. Parma}. $\textbf{7}$ ($2^{\ast}$) 1-99.

\bibitem{dophd} \textsc{Dolera, E.} (2010). Rapidity of convergence to equilibrium of the solution of the Boltzmann equation for Maxwellian molecules. Ph.D. thesis, Universit\`a degli Studi di Pavia.

\bibitem{doLom} \textsc{Dolera, E.} (2011). Spatially homogeneous Maxwellian molecules in a neighborhood of the equilibrium. In preparation.

\bibitem{dgr} \textsc{Dolera, E., Gabetta, E.} and \textsc{Regazzini, E.} (2009). Reaching the best possible rate of convergence to equilibrium for solutions of Kac's equation via central limit theorem. \emph{Ann. Appl. Probab}. $\textbf{19}$ 186-201.

\bibitem{dore} \textsc{Dolera, E.} and \textsc{Regazzini, E.} (2010). The role of the central limit theorem in discovering sharp rates of convergence to equilibrium for the solution of the Kac equation. \emph{Ann. Appl. Probab}. $\textbf{20}$ 430-461.

\bibitem{du} \textsc{Dudley, R.M.} (2002). \emph{Real Analysis and Probability}. Cambridge University Press, Cambridge.

\bibitem{elm} \textsc{Elmroth, T.} (1983). Global boundedness of moments of solutions of the Boltmann equation for forces of infinite range. \emph{Arch. Rational Mech. Anal}. $\textbf{82}$ 1-12.

\bibitem{fol} \textsc{Folland, G.B.} (1999). \emph{Real Analysis. Modern Techniques and Their Applications}. $2^{nd}$ ed. Wiley, NY.

\bibitem{flr} \textsc{Fortini, S., Ladelli, L.} and \textsc{Regazzini, E.} (1996). A central limit problem for partially exchangeable random variables. \emph{Theory Probab. Appl}. $\textbf{41}$ 224-246.

\bibitem{frgr} \textsc{Fristedt, B.} and \textsc{Grey, L.} (1997). \emph{A Modern Approach to Probability Theory}. Birkh\"{a}user, Boston.

\bibitem{gr6} \textsc{Gabetta, E.} and \textsc{Regazzini, E.} (2006). Some new results for McKean's graphs with applications to Kac's equation. \emph{J. Stat. Phys.} $\textbf{125}$ 947-974.

\bibitem{gr8} \textsc{Gabetta, E.} and \textsc{Regazzini, E.} (2008). Central limit theorem for the solution of the Kac equation. \emph{Ann. Appl. Probab.} $\textbf{18}$ 2320-2336.

\bibitem{hir} \textsc{Hirsch, M. W.} (1976). \emph{Differential Topology}. Springer-Verlag, NY.

\bibitem{ka} \textsc{Kallenberg, O.} (2002). \emph{Foundations of Modern Probability}. $2^{nd}$ ed. Springer-Verlag, NY.

\bibitem{max} \textsc{Maxwell, J.C.} (1867). On the dynamical theory of gases. \emph{Philos. Trans. Roy. Soc. London Ser.} \textbf{157} 49-88.

\bibitem{mck6} \textsc{McKean Jr, H. P.} (1966). Speed of approach to equilibrium for Kac's caricature of a Maxwellian gas. \emph{Arch. Rational Mech. Anal.} $\textbf{21}$ 343-367.

\bibitem{mck7} \textsc{McKean Jr, H. P.} (1967). An exponential formula for solving Boltzmann's equation for a Maxwellian gas. \emph{J. Combinatorial Theory} $\textbf{2}$ 358-382.

\bibitem{mor} \textsc{Morgenstern, D.} (1954). General existence and uniqueness proof for the spatially homogeneous solution of the Maxwell-Boltzmann equation in the case of Maxwellian molecules. \emph{Proc. Nat. Acad. Sci. USA} \textbf{40}
    719-721.

\bibitem{parth} \textsc{Parthasarathy, K.R.} (1967). \emph{Probability Measures on Metric Spaces}. Academic press, NY. Reprinted in 2005 by AMS Chelsea, Providence.

\bibitem{pe1} \textsc{Petrov, V. V.} (1995). \emph{Limit Theorems of Probability Theory. Sequences of Independent Random Variables}. Oxford University Press, NY.

\bibitem{puto} \textsc{Pulvirenti, A.} and \textsc{Toscani, G.} (1996). The theory of the nonlinear Boltzmann equation for Mawxell molecules in Fourier representation. \emph{Ann. Mat. Pura Appl.} $\textbf{171}$ 181-204.

\bibitem{stro} \textsc{Stroock, D.W.} (2011). \emph{Probability Theory, an Analytical View}. $2^{nd}$ ed. Cambridge University Press, Cambridge.

\bibitem{ta} \textsc{Tanaka, H.} (1978). Probabilistic treatment of the Boltzmann equation of Maxwellian molecules. \emph{Z. Wahrsch. Verw. Gebiete.} $\textbf{46}$ 67-105.

\bibitem{vil}  \textsc{Villani, C.} (2002). A review of mathematical topics in collisional kinetic theory. \emph{Handbook of Mathematical Fluid Dynamics}. (S. Friedlander and D. Serre eds). Vol. I, 71-305. North-Holland, Amsterdam.

\bibitem{wil} \textsc{Wild, E.} (1951). On Boltzmann's equation in the kinetic theory of gases. \emph{Proc. Cambridge Philos. Soc.} \textbf{47} 602-609.

\end{thebibliography}
\end{document}